\documentclass[11pt]{amsart}
\usepackage{amsmath,amssymb,amsfonts,url,mathptmx}
\numberwithin{equation}{section}

\newtheorem{thm}{Theorem}[section]
\newtheorem{lem}{Lemma}[section]
\newtheorem{rem}{Remark}[section]
\newtheorem{prop}{Proposition}[section]

\begin{document}
\title[Singular Liouville Equation]{Laplacian Vanishing Theorem for Quantized Singular Liouville Equation}

\author{Juncheng Wei}
\address{Department of Mathematics \\ Chinese University of Hong Kong, Shatin, NT, Hong Kong} \email{jcwei@math.ubc.ca }

\author{Lei Zhang} \footnote{The research of J. Wei is partially supported by NSERC of Canada. Lei Zhang is partially supported by a Simons Foundation Collaboration Grant}
\address{Department of Mathematics\\
        University of Florida\\
        1400 Stadium Rd\\
        Gainesville FL 32611}
\email{leizhang@ufl.edu}

\date{\today}

\begin{abstract} In this article we establish a vanishing theorem for singular Liouville equation with quantized singular source. If a blowup sequence tends to infinity near a quantized singular source and the blowup solutions violate the spherical Harnack inequality around the singular source (non-simple blow-ups), the Laplacian of a coefficient function must tend to zero. This seems to be the first second order estimates for Liouville equation with quantized sources and non-simple blow-ups.  This result as well as the key ideas of the proof would be extremely useful for various applications. 
\end{abstract}

\maketitle

\section{Introduction}
This is the third article in our series to study blowup solutions of
\begin{equation}\label{main-2}
\Delta u+|x|^{2N}\mathrm{H}(x)e^{u}=0, 
\end{equation}
in a neighborhood of the origin in $\mathbb R^2$. Here $\mathrm{H}$ is a positive smooth function and $N\in \mathbb N$ is a positive integer. Since the analysis is local in nature we focus the discussion in a neighborhood of the origin:
Let $\mathfrak{u}_k$ be a sequence of solutions of
\begin{equation}\label{t-u-k}
\Delta \mathfrak{u}_k(x)+|x|^{2N}\mathrm{H}_k(x)e^{\mathfrak{u}_k}=0, \quad \mbox{in}\quad B_{\tau}
\end{equation}
for some $\tau>0$ independent of $k$. $B_{\tau}$ is the ball centered at the origin with radius $\tau$.  In addition we postulate the usual assumptions on $\mathfrak{u}_k$ and $\mathrm{H}_k$:
For a positive constant $C$ independent of $k$, the following holds:
\begin{equation}\label{assumption-1}
\left\{\begin{array}{ll}
\|\mathrm{H}_k\|_{C^3(\bar B_{\tau})}\le C, \quad \frac 1C\le \mathrm{H}_k(x)\le C, \quad x\in \bar B_{\tau}, \\ \\
\int_{B_{\tau}} \mathrm{H}_k e^{\mathfrak{u}_k}\le C,\\  \\
|\mathfrak{u}_k(x)- \mathfrak{u}_k(y)|\le C, \quad \forall x,y\in \partial B_{\tau},
\end{array}
\right.
\end{equation}
and since we study the asymptotic behavior of blowup solutions around the singular source, we assume that there is no blowup point except at the origin:
\begin{equation}\label{assump-2}
\max_{K\subset\subset B_{\tau}\setminus \{0\}} \mathfrak{u}_k\le C(K).
\end{equation}

 If a sequence of solutions $\{u^k\}_{k=1}^{\infty}$ of (\ref{main-2}) satisfies
 $$\lim_{k\to \infty}u^k(x_k)=\infty,\quad \mbox{ for some $\bar x\in B_{\tau}$ and $x_k\to \bar x$,} $$
  we say $\{u^k\}$ is a sequence of bubbling solutions or blowup solutions, $\bar x$ is called a blowup point. The question we consider in this work is when $0$ is the only blowup point in a neighborhood of the origin, what vanishing theorems will the coefficient functions $\mathrm{H}_k$ satisfy?

One indispensable assumption is that the blowup solutions violate the spherical Harnack inequality around the origin:
\begin{equation}\label{no-sp-h}
\max_{x\in B_{\tau}} \mathfrak{u}_k(x)+2(1+N)\log |x|\to \infty,
\end{equation}
 It is also mentioned in literature ( see \cite{kuo-lin-jdg, wei-zhang-adv} ) that $0$ is called an non-simple blowup point. The main result of this article is

\begin{thm}\label{main-thm-1}
Let $\{\mathfrak{u}_k\}$ be a sequence of solutions of (\ref{t-u-k}) such that (\ref{assumption-1}),(\ref{assump-2}) hold and the spherical Harnack inequality is violated as in (\ref{no-sp-h}). Then along a sub-sequence
$$\lim_{k\to \infty}\Delta (\log \mathrm{H}_k)(0)=0. $$
\end{thm}

Theorem \ref{main-thm-1} is a continuation of our previous result in \cite{wei-zhang-plms}:

\medskip

\emph{Theorem A:  Let $\{\mathfrak{u}_k\}$ be a sequence of solutions of (\ref{t-u-k}) such that (\ref{assumption-1}),(\ref{assump-2}) and (\ref{no-sp-h}) hold. Then along a subsequence
$$\lim_{k\to \infty}\nabla (\log\mathrm{H}_k+\phi_k)(0)=0$$ where $\phi_k$ is defined as
\begin{equation}\label{phi-k}
\left\{\begin{array}{ll}
\Delta \phi_k(x)=0,\quad \mbox{in}\quad B_{\tau},\\ \\
\phi_k(x)=\mathfrak{u}_k(x)-\frac{1}{2\pi \tau}\int_{\partial B_{\tau}}\mathfrak{u}_kdS,\quad x\in \partial B_{\tau}.
\end{array}
\right.
\end{equation}
}
The equation (\ref{main-2}) comes from its equivalent form 
$$\Delta v+\mathrm{H}e^v=4\pi N\delta_0 $$
by using a logarithmic function to eliminate the Dirac mass on the right hand side. Since the strength of the Dirac mass is a multiple of $4\pi$, this type of singularity is called ``quantized". 
An equation with a quantized singular source is ubiquitous in literature. In particular
the following mean field equation defined on a Riemann surface $(M,g)$:
\begin{equation}\label{mean-sin}
\Delta_gu+\rho(\frac{h(x)e^{u(x)}}{\int_Mhe^{u}}-1)=4\pi\sum_{t=1}^M \alpha_t (\delta_{p_t}-1),
\end{equation}
 represents a conformal metric with prescribed conic singularities (see \cite{erem-3,tr-1,tr-2}), where $h$ is a positive smooth function, $\rho>0$ is a constant and the volume of $M$ is assumed to be $1$ for convenience, $\alpha_j>-1$ are constants as well. If the singular source is quantized, the equation is profoundly linked to Algebraic geometry, integrable system, number theory and complex Monge-Ampere equations (see \cite{chen-lin-last-cpam}). In Physics the main equation reveals key features of mean field limits of point vortices in the Euler flow \cite{caglioti-1,caglioti-2} and models in the Chern-Simons-Higgs theory \cite{jackiw} and in the electroweak theory \cite{ambjorn}, etc.

So far the non-simple bubbling situation has been observed in Liouville equation \cite{kuo-lin-jdg,bart3}, Liouville systems \cite{gu-zhang-1,gu-zhang-2, wu-zhang-siam} and fourth order equations \cite{ahmedou-wu-zhang}. The main theorem in this article would impact the study of these equations as well as some well known open questions in Monge-Ampere equation \cite{wei-zhang-adv}.

When compared with Theorem A, Theorem \ref{main-thm-1} is clearly more challenging in analysis. In fact the proof of Theorem A is a special case of one step of the proof of Theorem \ref{main-thm-1}. However, their major difference is on applications. Theorem \ref{main-thm-1} is significantly more influential for many reasons: First the main motivation to study equation (\ref{main-2}) is for equations or systems defined on manifolds. Usually blowup analysis near a singular point needs to reflect the curvature at the blowup point. In this respect Theorem \ref{main-thm-1} is directly related to the Gauss curvature at the blowup point. Second, the harmonic function in Theorem A causes inconvenience in application since it is generally hard to identify what the harmonic function is. On the other hand Theorem \ref{main-thm-1} is only involved with the Laplace of the coefficient function. This may lead to substantial advances in applications: In many degree counting programs one major difficulty is bubble-coalition, which means bubbling disks may collide into one point. The formation of bubbling disks tending to one point is accurately represented by equation (\ref{main-2}). Theorem \ref{main-thm-1} and its proof could be extremely useful to simplify blowup pictures.  Third, the proof of Theorem \ref{main-thm-1} is also important for proving uniqueness of bubbling solutions, and the results for Liouville equation with quantized singular sources is inspirational to many equations and systems with similar singular poles. Before our series of works most of the study of singular equations or systems focuses on non-quantized singular situations. However it is the ``quantized situations" that manifest profound connections to different fields of mathematics and Physics. Theorem \ref{main-thm-1} may be a starting point of multiple directions of exciting adventures. 

As a first application of Theorem \ref{main-thm-1} we present an advancement of the mean field equation (\ref{mean-sin}). 
Let $\Lambda$ be defined as
$$\Lambda=\{8\pi k+\sum_{j\in A}8\pi (1+\alpha_j);\quad k\in \mathbb N\cup \{0\}, \quad A\subset \{1,...,M\}\quad \} $$
where $\mathbb N$ is the set of natural numbers. 
Through the works of Bartoclucci-Tarantello \cite{bart-taran-jde-2,bart4}, Chen-Lin \cite{chen-lin-last-cpam} etc, an a priori estimate holds if $\rho\not \in \Lambda$. In other words, if $u^k$ is a sequence of blowup solutions with parameters $\rho^k$, the limit of $\rho^k$ is in $\Lambda$. 
Our second main theorem is
\begin{thm}\label{global-estimate}
Let $u_k$ be a sequence of blowup solutions of (\ref{mean-sin}) with parameters $\rho^k\to \rho\in \Lambda$, $h$ be a positive smooth function, $\alpha_1,...,\alpha_M>-1$ are constants. If at each quantized blowup point $p$ we have
\begin{equation}\label{curvature-a}
\Delta \log h(p)-2K(p)-4\pi\sum_{t=1}^M\alpha_t+\rho\neq 0, 
\end{equation}
where $K(p)$ is the Gauss curvature at $p$, all blowup points of $u_k$ are simple blowup points.
\end{thm}

The organization of the article is as follows. In section two we cite preliminary results related to the proof of the main theorem. Then in section three approximate the blowup solutions by a family of global solutions that agree with the blowup solutions at one local maximum point. This is crucial for our argument. Then we derive some intermediate estimates as preparation of more precise analysis. In section four we prove the first order estimates that cover the main result in \cite{wei-zhang-plms}. This section proves stronger result than \cite{wei-zhang-plms} and provides more detail. Finally in section five we take advantage of the result of the first order estimate and complete the proof of Theorem \ref{main-thm-1}. The proof of Theorem \ref{global-estimate} is placed in section six. The final section is an appendix that contains certain computations needed in the proof of Theorem \ref{main-thm-1}. 

{\bf Notation:} We will use $B(x_0,r)$ to denote a ball centered at $x_0$ with radius $r$. If $x_0$ is the origin we use $B_r$. $C$ represents a positive constant that may change from place to place.

\section{Preliminary discussions}
In the first stage of the proof of Theorem \ref{main-thm-1} we set up some notations and cite some preliminary results.
For simple notation we set
\begin{equation}\label{uk-d}
u_k(x)=\mathfrak{u}_k(x)-\phi_k(x), \quad \mbox{and}
\end{equation}
\begin{equation}\label{hk-d}
h_k(x)=\mathrm{H}_k(x)e^{\phi_k(x)}.
\end{equation}
to write the equation of $u_k$ as
\begin{equation}\label{eq-uk}
\Delta u_k(x)+|x|^{2N}h_k(x)e^{u_k}=0,\quad \mbox{ in }\quad B_{\tau}
\end{equation}
Without loss of generality we assume
\begin{equation}\label{rea-h}
\lim_{k\to \infty} h_k(0)=1.
\end{equation}

Obviously (\ref{no-sp-h}) is equivalent to
\begin{equation}\label{no-sp-h-u}
\max_{x\in B_{\tau}} u_k(x)+2(1+N)\log |x|\to \infty,
\end{equation}

It is well known \cite{kuo-lin-jdg, bart3} that $ u_k$ exhibits a non-simple blowup profile.  It is established in \cite{kuo-lin-jdg,bart3} that there are $N+1$ local maximum points of $ u_k$: $p_0^k$,....,$p_{N}^k$ and they are evenly distributed on $\mathbb S^1$ after scaling according to their magnitude: Suppose along a subsequence
$$\lim_{k\to \infty}p_0^k/|p_0^k|=e^{i\theta_0}, $$
then
$$\lim_{k\to \infty} \frac{p_l^k}{|p_0^k|}=e^{i(\theta_0+\frac{2\pi l}{N+1})}, \quad l=1,...,N. $$
For many reasons it is convenient to denote $|p_0^k|$ as $\delta_k$ and define $\mu_k$ as follows:
\begin{equation}\label{muk-dk}
\delta_k=|p_0^k|\quad \mbox{and }\quad \mu_k= u_k(p_0^k)+2(1+N)\log \delta_k.
\end{equation}
Also we use 
$$\epsilon_k=e^{-\frac 12 \mu_k} $$
to be the scaling factor most of the time. 
Since $p_l^k$'s are evenly distributed
around $\partial B_{\delta_k}$, standard results for Liouville equations around a regular blowup point can be applied to have $ u_k(p_l^k)= u_k(p_0^k)+o(1)$. Also, (\ref{no-sp-h}) gives $\mu_k\to \infty$. The interested readers may look into \cite{kuo-lin-jdg,bart3} for more detailed information.

Finally we shall use $E$ to denote a frequently appearing error term of the size $O(\delta_k^2)+O(\mu_ke^{-\mu_k})$.

\section{Approximating bubbling solutions by global solutions}

We write $p_0^k$ as $p_0^k=\delta_ke^{i\theta_k}$ and define $v_k$ as
\begin{equation}\label{v-k-d}
v_k(y)=u_k(\delta_k ye^{i\theta_k})+2(N+1)\log \delta_k,\quad |y|<\tau \delta_k^{-1}.
\end{equation}
If we write out each component, (\ref{v-k-d}) is
$$
v_k(y_1,y_2)=u_k(\delta_k(y_1\cos\theta_k-y_2\sin\theta_k),\delta_k(y_1\sin\theta_k+y_2\cos\theta_k))+2(1+N)\log \delta_k. $$
Then it is standard to verify that $v_k$ solves

\begin{equation}\label{e-f-vk}
\Delta v_k(y)+|y|^{2N}\mathfrak{h}_k(\delta_k y)e^{v_k(y)}=0,\quad |y|<\tau/\delta_k,
\end{equation}
where
\begin{equation}\label{frak-h}
\mathfrak{h}_k(x)=h_k(xe^{i\theta_k}),\quad |x|<\tau.
\end{equation}
Thus the image of $p_0^k$ after scaling is $Q_1^k=e_1=(1,0)$.
Let $Q_1^k$, $Q_2^k$,...,$Q_{N}^k$ be the images of $p_i^k$ $(i=1,...,N)$ after the scaling:
$$Q_l^k=\frac{p_l^k}{\delta_k}e^{-i\theta_k},\quad l=1,...,N. $$
 It is established by Kuo-Lin in \cite{kuo-lin-jdg} and independently by Bartolucci-Tarantello in \cite{bart3} that
\begin{equation}\label{limit-q}
\lim_{k\to \infty} Q_l^k=\lim_{k\to \infty}p_l^k/\delta_k=e^{\frac{2l\pi i}{N+1}},\quad l=0,....,N.
\end{equation}
Then it is proved in our previous work that ( see (3.13) in \cite{wei-zhang-adv})

$$Q_l^k-e^{\frac{2\pi l i}{N+1}}=O(\mu_ke^{-\mu_k})+O(|\nabla \log \mathfrak{h}_k(0)|\delta_k). $$
Using the rate of $\nabla \mathfrak{h}_k(0)$ in \cite{wei-zhang-adv} we have
\begin{equation}\label{Qm-close}
Q_l^k-e^{\frac{2\pi l i}{N+1}}=O(\mu_ke^{-\mu_k})+O(\delta_k^2).
\end{equation}
Choosing $3\epsilon>0$ small and independent of $k$, we can make disks centered at $Q_l^k$ with radius $3\epsilon$ (denoted as $B(Q_l^k,3\epsilon ) $) mutually disjoint. Let
\begin{equation}\label{v-muk}
\mu_k=\max_{B(Q_0^k,\epsilon)} v_k.
\end{equation}
Since $Q_l^k$ are evenly distributed around $\partial B_1$, it is easy to use standard estimates for single Liouville equations (\cite{zhangcmp,gluck,chenlin1}) to obtain
$$\max_{B(Q_l^k,\epsilon)}v_k=\mu_k+o(1),\quad l=1,...,N. $$

Let
\begin{equation}\label{def-Vk}
V_k(x)=\log \frac{e^{\mu_k}}{(1+\frac{e^{\mu_k}\mathfrak{h}_k(\delta_k e_1)}{8(1+N)^2}|y^{N+1}-e_1|^2)^2}.
\end{equation}
Clearly $V_k$ is a solution of
\begin{equation}\label{eq-for-Vk}
\Delta V_k+\mathfrak{h}_k(\delta_k e_1)|y|^{2N}e^{V_k}=0,\quad \mbox{in}\quad \mathbb R^2, \quad V_k(e_1)=\mu_k.
\end{equation}
This expression is based on the classification theorem of Prajapat-Tarantello \cite{prajapat}.

The estimate of $v_k(x)-V_k(x)$ is important for the main theorem of this article. For convenience we use
$$\beta_l=\frac{2\pi l}{N+1}, \quad \mbox{so}\,\, e_1=e^{i\beta_0}=Q_0^k,\quad
e^{i\beta_l}=Q_l^k+E,\,\,\mbox{ for }\,\, l=1,...,N. $$

\section{Vanishing of the first derivatives}
Our first goal is to prove the following vanishing rate for $\nabla \mathfrak{h}_k(0)$:
\begin{thm}\label{vanish-first-h}
\begin{equation}\label{vanish-first-tau}
\nabla (\log \mathfrak{h}_k)(0)=O(\delta_k\mu_k)
\end{equation}
\end{thm}

\noindent{\bf Proof of Theorem \ref{vanish-first-h}:}

Note that we have proved in \cite{wei-zhang-adv} that 
$$\nabla (\log \mathfrak{h}_k)(0)=O(\delta_k^{-1}\mu_ke^{-\mu_k})+O(\delta_k).$$
If $\delta_k\ge C\epsilon_k$, there is nothing to prove. So we assume that 
\begin{equation}\label{delta-small-1}
\delta_k=o(\epsilon_k). 
\end{equation}
By way of contradiction we assume that 
\begin{equation}\label{assum-delta-h}
|\nabla \mathfrak{h}_k(0)|/(\delta_k \mu_k)\to \infty. 
\end{equation}

Another observation is that based on (\ref{Qm-close}) we have 
$$\epsilon_k^{-1} | Q_l^k-e^{i\beta_l}|\le C\epsilon_k^{\epsilon},\quad l=0,...,N$$
for some small $\epsilon>0$. Thus $\xi_k$ tends to $U$ after scaling. We need this fact in our argument.

Under the assumption (\ref{delta-small-1}) we cite Proposition 3.1 of \cite{wei-zhang-plms}: 

\emph{Proposition 3.1 of \cite{wei-zhang-plms}: Let $l=0,...,N$ and $\delta$ be small so that $B(e^{i\beta_l},\delta)\cap B(e^{i\beta_s},\delta)=\emptyset$ for $l\neq s$.
In each $B(e^{i\beta_l},\delta)$
\begin{equation}\label{global-close}
|v_k(x)-V_k(x)|\le \left\{\begin{array}{ll}
C\mu_ke^{-\mu_k/2},\quad |x-e^{i\beta_l}|\le Ce^{-\mu_k/2}, \\
\\
C\frac{\mu_ke^{-\mu_k}}{|x-e^{i\beta_l}|}+O(\mu_k^2e^{-\mu_k}),\quad Ce^{-\mu_k/2}\le |x-e^{i\beta_l}|\le \delta.
\end{array}
\right.
\end{equation} }

\medskip

\begin{rem} We only need a re-scaled version of the Proposition above:  
\begin{equation}\label{vk-Vk-2}
|v_k(e^{i\beta_l}+\epsilon_ky)-V_k(e^{i\beta_l}+\epsilon_ky)|\le C\epsilon_k^{\epsilon} (1+|y|)^{-1},\quad 0<|y|<\tau \epsilon_k^{-1}.
\end{equation}
for some small constants $\epsilon>0$ and $\tau>0$ both independent of $k$,
\end{rem}

One major step in the proof of Theorem \ref{vanish-first-h} is the following estimate:

\begin{prop}\label{key-w8-8} Let $w_k=v_k-V_k$, then
$$|w_k(y)|\le C\tilde{\delta_k}, \quad y\in \Omega_k:=B(0,\tau \delta_k^{-1}), $$
where $\tilde{\delta_k}=|\nabla \mathfrak{h}_k(0)|\delta_k+\delta_k^2\mu_k$.
\end{prop}

\noindent{\bf Proof of Proposition \ref{key-w8-8}:}

Obviously we can assume that $|\nabla \mathfrak{h}_k(0)|\delta_k> 2\delta_k^2\mu_k$ because otherwise there is nothing to prove. 
Now we recall the equation for $v_k$ is (\ref{e-f-vk}),
$v_k$ is a constant on $\partial B(0,\tau \delta_k^{-1})$. Moreover $v_k(e_1)=\mu_k$. Recall that $V_k$ defined in (\ref{def-Vk}) satisfies
$$\Delta V_k+\mathfrak{h}_k(\delta_ke_1)|y|^{2N}e^{V_k}=0,\quad \mbox{in}\quad \mathbb R^2, \quad \int_{\mathbb R^2}|y|^{2N}e^{V_k}<\infty, $$
$V_k$ has its local maximums at $e^{i\beta_l}$ for $l=0,...,N$ and $V_k(e_1)=\mu_k$.
For $|y|\sim \delta_k^{-1}$,
$$V_k(y)=-\mu_k-4(N+1)\log \delta_k^{-1}+C+O(\delta_k^{N+1}). $$

Let $\Omega_k=B(0,\tau \delta_k^{-1})$,  we shall derive a precise, point-wise estimate of $w_k$ in $B_3\setminus \cup_{l=1}^{N}B(Q_l^k,\tau)$ where $\tau>0$ is a small number independent of $k$. Here we note that among $N+1$ local maximum points, we already have $e_1$ as a common local maximum point for both $v_k$ and $V_k$ and we shall prove that $w_k$ is very small in $B_3$ if we exclude all bubbling disks except the one around $e_1$. Before we carry out more specific computation we emphasize the importance of
\begin{equation}\label{control-e}
w_k(e_1)=|\nabla w_k(e_1)|=0.
\end{equation}
Now we write the equation of $w_k$ as
\begin{equation}\label{eq-wk}
\Delta w_k+\mathfrak{h}_k(\delta_k y)|y|^{2N}e^{\xi_k}w_k=(\mathfrak{h}_k(\delta_k e_1)-\mathfrak{h}_k(\delta_k y))|y|^{2N}e^{V_k}
\end{equation}
 in $\Omega_k$, where $\xi_k$ is obtained from the mean value theorem:
$$
e^{\xi_k(x)}=\left\{\begin{array}{ll}
\frac{e^{v_k(x)}-e^{V_k(x)}}{v_k(x)-V_k(x)},\quad \mbox{if}\quad v_k(x)\neq V_k(x),\\
\\
e^{V_k(x)},\quad \mbox{if}\quad v_k(x)=V_k(x).
\end{array}
\right.
$$
An equivalent form is
\begin{equation}\label{xi-k}
e^{\xi_k(x)}=\int_0^1\frac d{dt}e^{tv_k(x)+(1-t)V_k(x)}dt=e^{V_k(x)}\big (1+\frac 12w_k(x)+O(w_k(x)^2)\big ).
\end{equation}
For convenience we write the equation for $w_k$ as
\begin{equation}\label{eq-wk-2}
\Delta w_k+\mathfrak{h}_k(\delta_k y)|y|^{2N}e^{\xi_k}w_k=\delta_k\nabla \mathfrak{h}_k(\delta_k e_1)\cdot (e_1-y)|y|^{2N}e^{V_k}+E_1
\end{equation}
where $$E_1=O(\delta_k^2)|y-e_1|^2|y|^{2N}e^{V_k},\quad y\in \Omega_k. $$
Note that the oscillation of $w_k$ on $\partial \Omega_k$ is $O(\delta_k^{N+1})$, which all comes from the oscillation of $V_k$. 

Let $M_k=\max_{x\in \bar \Omega_k}|w_k(x)|$. We shall get a contradiction by assuming $M_k/\tilde{\delta_k}\to \infty$. This assumption implies
\begin{equation}\label{big-mk}
M_k/(\delta_k^2 \mu_k)\to \infty. 
\end{equation}
Set
$$\tilde w_k(y)=w_k(y)/M_k,\quad x\in \Omega_k. $$
Clearly $\max_{x\in \Omega_k}|\tilde w_k(x)|=1$. The equation for $\tilde w_k$ is
\begin{equation}\label{t-wk}
\Delta \tilde w_k(y)+|y|^{2N}\mathfrak{h}_k(\delta_k e_1)e^{\xi_k}\tilde w_k(y)=a_k\cdot (e_1-y)|y|^{2N}e^{V_k}+\tilde E_1,
\end{equation}
in $\Omega_k$,
where $a_k=\delta_k\nabla \mathfrak{h}_k(0)/M_k\to 0$,
\begin{equation}\label{t-ek}
\tilde E_1=o(1)|y-e_1|^2|y|^{2N}e^{V_k},\quad y\in \Omega_k.
\end{equation}
Also on the boundary, since $M_k/\tilde \delta_k\to \infty$, we have 
\begin{equation}\label{w-bry}
\tilde w_k=C+o(1/\mu_k),\quad \mbox{on}\quad \partial \Omega_k. 
\end{equation}

 By Proposition 3.1 of \cite{wei-zhang-plms}
\begin{equation}\label{xi-V-c}
\xi_k(e_1+\epsilon_k z)=V_k(e_1+\epsilon_kz)+O(\epsilon_k^{\epsilon})(1+|z|)^{-1}
\end{equation}

Since $V_k$ is not exactly symmetric around $e_1$, we shall replace the re-scaled version of $V_k$ around $e_1$ by a radial function.
Let $U_k$ be solutions of
\begin{equation}\label{global-to-use}
\Delta U_k+\mathfrak{h}_k(\delta_ke_1)e^{U_k}=0,\quad \mbox{in}\quad \mathbb R^2, \quad U_k(0)=\max_{\mathbb R^2}U_k=0.
\end{equation}
By the classification theorem of Caffarelli-Gidas-Spruck \cite{CGS} we have
$$U_k(z)=\log \frac{1}{(1+\frac{\mathfrak{h}_k(\delta_ke_1)}{8}|z|^2)^2}$$
and standard refined estimates yield (see \cite{chenlin1,zhangcmp,gluck})
\begin{equation}\label{Vk-rad}
V_k(e_1+\epsilon_k z)+2\log \epsilon_k=U_k(z)+O(\epsilon_k)|z|+O(\mu_k^2\epsilon_k^2).
\end{equation}
Also we observe that
\begin{equation}\label{log-rad}
\log |e_1+\epsilon_k z|=O(\epsilon_k)|z|.
\end{equation}

Thus, the combination of (\ref{xi-V-c}), (\ref{Vk-rad}) and (\ref{log-rad}) gives
\begin{align}\label{xi-U}
&2N\log |e_1+\epsilon_kz|+\xi_k(e_1+\epsilon_k z)+2\log \epsilon_k-U_k(z)\\
=&O(\epsilon_k^{\epsilon})(1+|z|)\quad 0\le |z|<\delta_0 \epsilon_k^{-1}.
 \nonumber
\end{align}
for a small $\epsilon>0$ independent of $k$. 
Since we shall use the re-scaled version, based on (\ref{xi-U}) we have
\begin{equation}\label{xi-eU}
\epsilon_k^2 |e_1+\epsilon_k z|^{2N}e^{\xi_k(e_1+\epsilon_k z)}
= e^{U_k(z)}+O( \epsilon_k^{\epsilon})(1+|z|)^{-3}
\end{equation}
Here we note that the estimate in (\ref{xi-U}) is not optimal.  In the following we shall put the proof of Proposition \ref{key-w8-8} into a few estimates. In the first estimate we prove

\begin{lem}\label{w-around-e1} For $\delta>0$ small and independent of $k$,
\begin{equation}\label{key-step-1}
\tilde w_k(y)=o(1),\quad \nabla \tilde w_k=o(1) \quad \mbox{in}\quad B(e_1,\delta)\setminus B(e_1,\delta/8)
\end{equation}
where $B(e_1,3\delta)$ does not include other blowup points.
\end{lem}

\noindent{\bf Proof of Lemma  \ref{w-around-e1}:}

If (\ref{key-step-1}) is not true, we have, without loss of generality that $\tilde w_k\to c>0$. This is based on the fact that $\tilde w_k$ tends to a global harmonic function with removable singularity. So $\tilde w_k$ tends to constant. Here we assume $c>0$ but the argument for $c<0$ is the same. Let
\begin{equation}\label{w-ar-e1}
W_k(z)=\tilde w_k(e_1+\epsilon_kz), \quad \epsilon_k=e^{-\frac 12 \mu_k},
\end{equation}
then if we use $W$ to denote the limit of $W_k$, we have
$$\Delta W+e^UW=0, \quad \mathbb R^2, \quad |W|\le 1, $$
and $U$ is a solution of $\Delta U+e^U=0$ in $\mathbb R^2$ with $\int_{\mathbb R^2}e^U<\infty$. Since $0$ is the local maximum of $U$,
$$U(z)=\log \frac{1}{(1+\frac 18|z|^2)^2}. $$
Here we further claim that $W\equiv 0$ in $\mathbb R^2$ because $W(0)=|\nabla W(0)|=0$, a fact well known based on the classification of the kernel of the linearized operator. Going back to $W_k$, we have
$$W_k(z)=o(1),\quad |z|\le R_k \mbox{ for some } \quad R_k\to \infty. $$

Based on the expression of $\tilde w_k$, (\ref{Vk-rad}) and (\ref{xi-eU}) we write the equation of $W_k$ as
\begin{equation}\label{e-Wk}
\Delta W_k(z)+\mathfrak{h}_k(\delta_ke_1)e^{U_k(z)}W_k(z)=E_2^k,
\end{equation}
for $|z|<\delta_0 \epsilon_k^{-1}$ where a crude estimate of the error term $E_2^k$ is
\begin{equation*}
E_2^k(z)=o(1)\epsilon_k^{\epsilon}(1+|z|)^{-3}.
\end{equation*}

Let
\begin{equation}\label{for-g0}
g_0^k(r)=\frac 1{2\pi}\int_0^{2\pi}W_k(r,\theta)d\theta.
\end{equation}
Then clearly $g_0^k(r)\to c>0$ for $r\sim \epsilon_k^{-1}$.
 The equation for $g_0^k$ is
\begin{align*}
&\frac{d^2}{dr^2}g_0^k(r)+\frac 1r \frac{d}{dr}g_0^k(r)+\mathfrak{h}_k(\delta_ke_1)e^{U_k(r)}g_0^k(r)=\tilde E_0^k(r)\\
&g_0^k(0)=\frac{d}{dr}g_0^k(0)=0.
\end{align*}
where $\tilde E_0^k(r)$ has the same upper bound as that of $E_2^k(r)$:
$$|\tilde E_0^k(r)|\le o(1)\epsilon_k^{\epsilon}(1+r)^{-3}. $$

For the homogeneous equation, the two fundamental solutions are known: $g_{01}$, $g_{02}$, where
$$g_{01}=\frac{1-c_1r^2}{1+c_1r^2},\quad c_1=\frac{\mathfrak{h}_k(\delta_ke_1)}8.$$
By the standard reduction of order process, $g_{02}(r)=O(\log r)$ for $r>1$.
Then it is easy to obtain, assuming $|W_k(z)|\le 1$, that
\begin{align*}
|g_0(r)|\le C|g_{01}(r)|\int_0^r s|\tilde E_0^k(s) g_{02}(s)|ds+C|g_{02}(r)|\int_0^r s|g_{01}(s)\tilde E_0^k(s)|ds\\
\le C\epsilon_k^{\epsilon}\log (2+r). \quad 0<r<\delta_0 \epsilon_k^{-1}.
\end{align*}
Clearly this is a contradiction to (\ref{for-g0}). We have proved $c=0$, which means $\tilde w_k=o(1)$ in $B(e_1, \delta_0)\setminus B(e_1, \delta_0/8)$.
Then it is easy to use the equation for $\tilde w_k$ and standard Harnack inequality to prove
$\nabla \tilde w_k=o(1)$ in the same region.
Lemma \ref{w-around-e1} is established. $\Box$

\medskip

The second estimate is a more precise description of $\tilde w_k$ around $e_1$:
\begin{lem}\label{t-w-1-better} For any given $\sigma\in (0,1)$ there exists $C>0$ such that
\begin{equation}\label{for-lambda-k}
|\tilde w_k(e_1+\epsilon_kz)|\le C\epsilon_k^{\sigma} (1+|z|)^{\sigma},\quad 0<|z|<\tau \epsilon_k^{-1}.
\end{equation}
for some $\tau>0$.
\end{lem}

\begin{rem} Lemma \ref{t-w-1-better} is an intermediate estimate for $\tilde w_k$. Eventually we shall improve (\ref{for-lambda-k}) to an error with the leading term $o(\epsilon_k)$.
\end{rem}

\noindent{\bf Proof of Lemma \ref{t-w-1-better}:} Let $W_k$ be defined as in (\ref{w-ar-e1}). In order to obtain a better estimate we need to write the equation of $W_k$ more precisely than (\ref{e-Wk}):
\begin{equation}\label{w-more}
\Delta W_k+\mathfrak{h}_k(\delta_ke_1)e^{\Theta_k}W_k=E_3^k(z), \quad z\in \Omega_{Wk}
\end{equation}
where
$\Theta_k$ is defined by
$$e^{\Theta_k(z)}=|e_1+\epsilon_k z|^{2N}e^{\xi_k(e_1+\epsilon_kz)+2\log \epsilon_k}, $$
$\Omega_{Wk}=B(0,\tau \epsilon_k^{-1})$ and $E_3^k(z)$ satisfies
$$E_3^k(z)=O(\epsilon_k)(1+|z|)^{-3},\quad z\in \Omega_{Wk}. $$
Here we observe that by Lemma \ref{w-around-e1}  $W_k=o(1)$ on $\partial \Omega_{Wk}$. 
Let
$$\Lambda_k=\max_{z\in \Omega_{Wk}}\frac{|W_k(z)|}{\epsilon_k^{\sigma}(1+|z|)^{\sigma}}. $$
If (\ref{for-lambda-k}) does not hold, $\Lambda_k\to \infty$ and we use $z_k$ to denote where $\Lambda_k$ is attained. Note that because of the smallness of $W_k$ on $\partial \Omega_{Wk}$, $z_k$ is an interior point. Let
$$g_k(z)=\frac{W_k(z)}{\Lambda_k (1+|z_k|)^{\sigma}\epsilon_k^{\sigma}},\quad z\in \Omega_{Wk}, $$
we see immediately that
\begin{equation}\label{g-sub-linear}
|g_k(z)|=\frac{|W_k(z)|}{\epsilon_k^{\sigma}\Lambda_k(1+|z|)^{\sigma}}\cdot \frac{(1+|z|)^{\sigma}}{(1+|z_k|)^{\sigma}}\le  \frac{(1+|z|)^{\sigma}}{(1+|z_k|)^{\sigma}}.
\end{equation}
Note that $\sigma$ can be as close to $1$ as needed. The equation of $g_k$ is
$$\Delta g_k(z)+\mathfrak{h}_k(\delta_k e_1)e^{\Theta_k}g_k=o(\epsilon_k^{1-\sigma})\frac{(1+|z|)^{-3}}{(1+|z_k|)^{\sigma}}, \quad \mbox{in}\quad \Omega_{Wk}. $$
Then we can obtain a contradiction to $|g_k(z_k)|=1$ as follows: If $\lim_{k\to \infty}z_k=P\in \mathbb R^2$, this is not possible because that fact that $g_k(0)=|\nabla g_k(0)|=0$ and the sub-linear growth of $g_k$ in (\ref{g-sub-linear}) implies that $g_k\to 0$ over any compact subset of $\mathbb R^2$ (see \cite{chenlin1,zhangcmp}). So we have $|z_k|\to \infty$. But this would lead to a contradiction again by using the Green's representation of $g_k$:
\begin{align} \label{temp-1}
&\pm 1=g_k(z_k)=g_k(z_k)-g_k(0)\\
&=\int_{\Omega_{k,1}}(G_k(z_k,\eta)-G_k(0,\eta))(\mathfrak{h}_k(\delta_k e_1)e^{\Theta_k}g_k(\eta)+o(\epsilon_k^{1-\sigma})\frac{(1+|\eta |)^{-3}}{(1+|z_k|)^{\sigma}})d\eta+o(1).\nonumber
\end{align}
where $G_k(y,\eta)$ is the Green's function on $\Omega_{Wk}$ and $o(1)$ in the equation above comes from the smallness of $W_k$ on $\partial \Omega_{Wk}$. Let $L_k=\tau\epsilon_k^{-1}$, the expression of $G_k$ is 
$$G_k(y,\eta)=-\frac{1}{2\pi}\log |y-\eta|+\frac 1{2\pi}\log (\frac{|\eta |}{L_k}|\frac{L_k^2\eta}{|\eta |^2}-y|). $$
$$G_k(z_k,\eta)-G_k(0,\eta)=-\frac{1}{2\pi}\log |z_k-\eta |+\frac 1{2\pi}\log |\frac{z_k}{|z_k|}-\frac{\eta z_k}{L_k^2}|+\frac 1{2\pi}\log |\eta |. $$
Using this expression in (\ref{temp-1}) we obtain from elementary computation that the right hand side of (\ref{temp-1}) is $o(1)$, a contradiction to $|g_k(z_k)|=1$. Lemma \ref{t-w-1-better} is
established. $\Box$

\medskip

The smallness of $\tilde w_k$ around $e_1$ can be used to obtain the following third key estimate:
\begin{lem}\label{small-other}
\begin{equation}\label{key-step-2}
\tilde w_k=o(1)\quad \mbox{in}\quad B(e^{i\beta_l},\tau)\quad l=1,..,N.
\end{equation}
\end{lem}

\noindent{\bf Proof of Lemma \ref{small-other}:}
We abuse the notation $W_k$ by defining it as
$$W_k(z)=\tilde w_k(e^{i\beta_l}+\epsilon_k z),\quad z\in \Omega_{k,l}:=B(0,\tau \epsilon_k^{-1}). $$
Here we point out that based on (\ref{Qm-close}) and (\ref{delta-small-1}) we have $\epsilon_k^{-1}|Q_l^k-e^{i\beta_l}|\to 0$. So the scaling around $e^{i\beta_l}$ or $Q_l^k$ does not affect the
limit function.

$$\epsilon_k^2 |e^{i\beta_l}+\epsilon_kz|^{2N}\mathfrak{h}_k(\delta_ke_1)e^{\xi_k(e^{i\beta_l}+\epsilon_kz)}\to e^{U(z)} $$
where $U(z)$ is a solution of
$$\Delta U+e^U=0,\quad \mbox{in}\quad \mathbb R^2, \quad \int_{\mathbb R^2}e^U<\infty. $$
Here we recall that $\lim_{k\to \infty} \mathfrak{h}_k(\delta_k e_1)=1$.
Since $W_k$ converges to a solution of the linearized equation:
$$\Delta W+e^UW=0, \quad \mbox{in}\quad \mathbb R^2. $$
$W$ can be written as a linear combination of three functions:
$$W(x)=c_0\phi_0+c_1\phi_1+c_2\phi_2, $$
where
$$\phi_0=\frac{1-\frac 18 |x|^2}{1+\frac 18 |x|^2} $$
$$\phi_1=\frac{x_1}{1+\frac 18 |x|^2},\quad \phi_2=\frac{x_2}{1+\frac 18|x|^2}. $$
The remaining part of the proof consisting of proving $c_0=0$ and $c_1=c_2=0$. First we prove $c_0=0$.

\noindent{\bf Step one: $c_0=0$.}
First we write the equation for $W_k$ in a convenient form. Since
$$|e^{i\beta_l}+\epsilon_kz|^{2N}\mathfrak{h}_k(\delta_ke_1)=\mathfrak{h}_k(\delta_ke_1)+O(\epsilon_k z),$$
and
$$\epsilon_k^2e^{\xi_k(e^{i\beta_l}+\epsilon_kz)}=e^{U_k(z)}+O(\epsilon_k^{\epsilon})(1+|z|)^{-3}. $$
Based on (\ref{t-wk}) we write the equation for $W_k$ as
\begin{equation}\label{around-l}
\Delta W_k(z)+\mathfrak{h}_k(\delta_ke_1)e^{U_k}W_k=E_l^k(z)
\end{equation}
where
$$E_l^k(z)=O(\epsilon_k^{\epsilon})(1+|z|)^{-3}\quad \mbox{in}\quad \Omega_{k,l}.$$
In order to prove $c_0=0$, the key is to control the derivative of $W_0^k(r)$ where
$$W_0^k(r)=\frac 1{2\pi r}\int_{\partial B_r} W_k(re^{i\theta})dS, \quad 0<r<\tau \epsilon_k^{-1}. $$
To obtain a control of $\frac{d}{dr}W_0^k(r)$ we use $\phi_0^k(r)$ as the radial solution of
$$\Delta \phi_0^k+\mathfrak{h}_k(\delta_k e_1)e^{U_k}\phi_0^k=0, \quad \mbox{in }\quad \mathbb R^2. $$

When $k\to \infty$, $\phi_0^k\to c_0\phi_0$. Thus using the equation for $\phi_0^k$ and $W_k$, we have
\begin{equation}\label{c-0-pf}
\int_{\partial B_r}(\partial_{\nu}W_k\phi_0^k-\partial_{\nu}\phi_0^kW_k)=o(\epsilon_k^{\epsilon}). \end{equation}

Thus from (\ref{c-0-pf}) we have
\begin{equation}\label{W-0-d}
\frac{d}{dr}W_0^k(r)=\frac{1}{2\pi r}\int_{\partial B_r}\partial_{\nu}W_k=o(\epsilon_k^{\epsilon})/r+O(1/r^3),\quad 1<r<\tau \epsilon_k^{-1}.
\end{equation}
Since we have known that
$$W_0^k(\tau \epsilon_k^{-1})=o(1). $$
By the fundamental theorem of calculus we have
$$W_0^k(r)=W_0^k(\tau\epsilon_k^{-1})+\int_{\tau \epsilon_k^{-1}}^r(\frac{o(\epsilon_k^{\epsilon})}{s}+O(s^{-3}))ds=O(1/r^2)+O(\epsilon_k^{\epsilon}\log
\frac{1}{\epsilon_k}) $$
for $r\ge 1$. Thus
$c_0=0$ because $W_0^k(r)\to c_0\phi_0$, which means when $r$ is large, it is $-c_0+O(1/r^2)$.

\medskip

\noindent{\bf Step two $c_1=c_2=0$}.
 We first observe that Lemma \ref{small-other} follows from this. Indeed, once we have proved $c_1=c_2=c_0=0$ around each $e^{i\beta_l}$, it is easy to use maximum principle to prove $\tilde w_k=o(1)$ in $B_3$ using $\tilde w_k=o(1)$ on $\partial B_3$ and the Green's representation of $\tilde w_k$. The smallness of $\tilde w_k$ immediately implies $\tilde w_k=o(1)$ in $B_R$ for any fixed $R>>1$. Outside $B_R$, a crude estimate of $v_k$ is
  $$v_k(y)\le -\mu_k-4(N+1)\log |y|+C, \quad 3<|y|<\tau \delta_k^{-1}. $$
  Using this and the Green's representation of $w_k$ we can first observe that the oscillation on each $\partial B_r$ is $o(1)$ ($R<r<\tau \delta_k^{-1}/2$) and then by the Green's representation of $\tilde w_k$ and fast decay rate of $e^{V_k}$ we obtain $\tilde w_k=o(1)$ in $\overline{B(0,\tau \delta_k^{-1})}$. A contradiction to $\max |\tilde w_k|=1$.

 There are $N+1$ local maximums with one of them being $e_1$. Correspondingly there are $N+1$ global solutions $V_{l,k}$ that
approximate $v_k$ accurately near $Q_l^k$ for $l=0,...,N$. Note that $Q_0^k=e_1$. For $V_{l,k}$ the expression is
$$V_{l,k}=\log \frac{e^{\mu_l^k}}{(1+\frac{e^{\mu_l^k}}{D_l^k}|y^{N+1}-(e_1+p_l^k)|^2)^2},\quad l=0,...,N, $$
where $p_l^k=E$ and
\begin{equation}\label{def-D}
D_l^k=8(N+1)^2/\mathfrak{h}_k(\delta_kQ_l^k).
\end{equation}
The equation that $V_{l,k}$ satisfies is 
$$\Delta V_{l,k}+|y|^{2N}\mathfrak{h}_k(\delta_k Q_l^k)e^{V_{l,k}}=0,\quad \mbox{in}\quad \mathbb R^2. $$
Since $v_k$ and $V_{l,k}$ have the same common local maximum at $Q_l^k$, it is easy to see that
\begin{equation}\label{ql-exp}
Q_l^k=e^{i\beta_l}+\frac{p_l^ke^{i\beta_l}}{N+1}+O(|p_l^k|^2),\quad \beta_l=\frac{2l\pi}{N+2}.
\end{equation}
Let $M_{l,k}$ be the maximum of $|v_k-V_{l,k}|$ and we claim that all these $M_{l,k}$ are comparable:
\begin{equation}\label{M-comp}
M_{l,k}\sim M_{s,k},\quad \forall s\neq l.
\end{equation}
The proof of (\ref{M-comp}) is as follows: We use $L_{s,l}$ to denote the limit of $(v_k-V_{l,k})/M_{l,k}$ around $Q_s^k$:
$$\frac{(v_k-V_{l,k})(Q_s^k+\epsilon_kz)}{M_{l,k}}=L_{s,l}+o(1),\quad |z|\le \tau \epsilon_k^{-1} $$
where
$$ L_{s,l}=c_{1,s,l}\frac{z_1}{1+\frac 18 |z|^2}+c_{2,s,l}\frac{z_2}{1+\frac 18 |z|^2},\quad \mbox{and}\quad L_{l,l}=0, \quad s=0,...,N. $$
If all $c_{1,s,l}$ and $c_{2,s,l}$ are zero for a fixed $l$, we can obtain a contradiction just like the beginning of step two. So at least one of them is not zero.
For each $s\neq l$, by Lemma \ref{t-w-1-better} we have
\begin{equation}\label{Q-bad}
v_k(Q_s^k+\epsilon_kz)-V_{s,k}(Q_s^k+\epsilon_kz)=O(\epsilon_k^{\sigma})(1+|z|)^{\sigma} M_{s,k},\quad |z|<\tau \epsilon_k^{-1}.
\end{equation}
Let $M_k=\max_{i}M_{i,k}$ ($i=0,...,N$) and we suppose $M_k=M_{l,k}$. Then to determine $L_{s,l}$ we see that
\begin{align*}
    &\frac{v_k(Q_s^k+\epsilon_k z)-V_{l,k}(Q_s^k+\epsilon_kz)}{M_k}\\
    =&o(\epsilon_k^{\sigma})(1+|z|)^{\sigma}+\frac{V_{s,k}(Q_s^k+\epsilon_kz)-V_{l,k}(Q_s^k+\epsilon_kz)}{M_k}. 
\end{align*}
This expression says that $L_{s,l}$ is mainly determined by the difference of two global solutions $V_{s,k}$ and $V_{l,k}$. In order to obtain a contradiction to our assumption we will put the difference in several terms. The main idea in this part of the reasoning is that ``first order terms" tell us what the kernel functions should be, then the ``second order terms" tell us where the pathology is. 

We write $V_{s,k}(y)-V_{l,k}(y)$ as
$$V_{s,k}(y)-V_{l,k}(y)=\mu_s^k-\mu_l^k+2A-A^2+O(|A|^3) $$
where
$$A(y)=\frac{\frac{e^{\mu_l^k}}{D_l^k}|y^{N+1}-e_1-p_l^k|^2-\frac{e^{\mu_s^k}}{D_s^k}|y^{N+1}-e_1-p_s^k|^2}{1+\frac{e^{\mu_s^k}}{D_s^k}|y^{N+1}-e_1-p_s^k|^2}.$$
Here for convenience we abuse the notation $\epsilon_k$ by assuming $\epsilon_k=e^{-\mu_s^k/2}$. Note that $\epsilon_k=e^{-\mu_t^k/2}$ for some $t$, but it does not matter which $t$ it is. From $A$ we claim that 
\begin{align}\label{late-1}
&V_{s,k}(Q_s^k+\epsilon_kz)-V_{l,k}(Q_s^k+\epsilon_kz)\\
=&\phi_1+\phi_2+\phi_3+\phi_4+\mathfrak{R},\nonumber
\end{align}
where
\begin{align*}
&\phi_1=(\mu_s^k-\mu_l^k)(1-\frac{(N+1)^2}{D_s^k}|z+O(\epsilon_k)|z|^2|^2)/B, \\
&\phi_2=\frac{2(N+1)^2}{D_s^k} \delta_k\nabla\mathfrak{h}_k(\delta_k Q_s^k)(Q_l^k-Q_s^k)|z|^2/B\\
&\phi_3=\frac{4(N+1)}{D_s^k B}Re((z+O(\epsilon_k|z|^2))(\frac{\bar p_s^k-\bar p_l^k}{\epsilon_k}e^{-i\beta_s}))
\\
&\phi_4=\frac{|p_s^k-p_l^k|^2}{\epsilon_k^2}\bigg (\frac{2}{D_s^k  B}-\frac{2(N+1)^2|z|^2}{D_s^2 B^2}
-\frac{2(N+1)^2}{D_s^2B^2}|z|^2\cos (2\theta-2\theta_{st}-2\beta_s) \bigg ),\\
&B=1+\frac{(N+1)^2}{D_s^k}|z+O(\epsilon_k|z|^2)|^2,
\end{align*}
and $\mathfrak{R}_k$ is the collections of other insignificant terms.  Here we briefly explain the roles of each term. $\phi_1$ corresponds to the radial solution in the kernel of the linearized operator of the global equation. In other words, $\phi_1^k/M_k$ should tend to zero because in step one we have proved $c_0=0$. $\phi_2^k/M_k$ is the combination of the two other functions in the kernel. $\phi_4$ is the second order term which will play a leading role later. $\phi_3^k$ comes from the difference of $\mathfrak{h}_k$ at $Q_l^k$ and $Q_s^k$. The derivation of (\ref{late-1}) is as follows: First by the expression of $Q_s^k$ in (\ref{ql-exp}) we have 
$$y^{N+1}=1+p_s^k+(N+1)\epsilon_kze^{-i\beta_s}+O(\epsilon_k^2)|z|^2, $$
where $y=Q_s^k+\epsilon_k z$. Then
$$|y^{N+1}-e_1-p_s^k|^2=(N+1)^2\epsilon_k^2|z+O(\epsilon_k)z^2|^2+O(\epsilon_k^3)|z|^3 $$
$$|y^{N+1}-e_1-p_l^k|^2=(N+1)^2\epsilon_k^2|z+\frac{(p_s^k-p_l^k)e^{i\beta_s}}{(N+1)\epsilon_k }+O(\epsilon_k)|z|^2|^2+O(\epsilon_k^3|z|^3). $$
Next by the definition of $D_s^k$ in (\ref{def-D}) 
$$\frac{D_s^k-D_l^k}{D_l^k}=\delta_k \nabla (\log \mathfrak{h}_k)(0)\cdot (Q_l^k-Q_s^k)+O(\delta_k^2).$$
\begin{align}\label{tem-5}
\frac{e^{\mu_l^k-\mu_s^k}}{D_l^k}&=\frac{1}{D_s^k}(1+\frac{D_s^k-D_l^k}{D_l^k}+\mu_l^k-\mu_s^k+O(\mu_l^k-\mu_s^k)^2+O(\delta_k^2)). \nonumber\\
&=\frac 1{D_s^k}(1+\delta_k\nabla\log \mathfrak{h}_k(0)\cdot (Q_l^k-Q_s^k)+\mu_l^k-\mu_s^k+O(\mu_l^k-\mu_s^k)^2+O(\delta_k^2)).
\end{align}
Then the expression of $A$ is ( for simplicity we omit $k$ in some notations)
\begin{align*}
    A=\bigg (\frac{e^{\mu_l-\mu_s}}{D_l}(N+1)^2\big (|z|^2+2Re(z\frac{\overline{p_s-p_l}}{\epsilon_k(N+1)}e^{-i\beta_s})+\frac{|p_s-p_l|^2}{(N+1)^2\epsilon_k^2}+O(\epsilon_k|z|^3)\big )\\
    -\frac{(N+1)^2}{D_s}(|z|^2+O(\epsilon_k)|z|^3)\bigg )/B
\end{align*}
After using (\ref{tem-5}) we have 
\begin{align}\label{exp-A-2}
A=\bigg (\frac{1}{D_s}(\delta_k\nabla(\log \mathfrak{h}_k(0)(Q_l-Q_s)+\mu_l-\mu_s+O(\mu_l-\mu_s)^2)(N+1)^2|z|^2\\
+2Re(z\frac{\bar p_s-\bar p_l}{\epsilon_k}e^{-i\beta_s})(N+1)\frac 1{D_s}+
\frac{|p_s-p_l|^2}{\epsilon_k^2D_s}+O(\epsilon_k)|z|^3+O(\delta_k^2)|z|^2\bigg )/B. \nonumber
\end{align}

\begin{equation}\label{exp-a-squ}
A^2=
(\frac{(N+1)^2}{D_s^2}4(Re(z\frac{\bar p_s-\bar p_l}{\epsilon_k}e^{-i\beta_s})^2)/B^2
+\mbox{other terms}.
\end{equation}
The numerator of $A^2$ has the following leading term:
$$\frac{(N+1)^2}{D_s^2}\bigg (2|z|^2(\frac{|p_s-p_l|}{\epsilon_k})^2\big (1+2\cos(2\theta-2\theta_{sl})\big ) \bigg )$$
where $z=|z|e^{i\theta}$, $p_s-p_l=|p_s-p_l|e^{i\theta_{sl}}$.
Using these expressions we can obtain (\ref{late-1}) by direct computation. 
Here $\phi_1$, $\phi_3$ correspond to solutions to the linearized operator. Here we note that if we set $\epsilon_{l,k}=e^{-\mu_l^k/2}$, there is no essential difference between $\epsilon_{l,k}$ and $\epsilon_k=e^{-\frac 12\mu_{1,k}}$ because $\epsilon_{l,k}=\epsilon_k+O(\epsilon_k E)$. If $|\mu_{s,k}-\mu_{l,k}|/M_k\ge C$ there is no way to obtain a limit in the form of $L_{s,l}$ mentioned before. Thus we must have $|\mu_{s,k}-\mu_{l,k}|/M_k\to 0$. After simplification (see $\phi_3$ of (\ref{late-1})) we have
\begin{align}\label{c-12}
c_{1,s,l}=\lim_{k\to \infty}\frac{|p_s^k-p_l^k|}{2(N+1)M_k\epsilon_k}\cos(\beta_s+\theta_{sl}),\\
c_{2,s,l}=\lim_{k\to \infty}
\frac{|p_s^k-p_l^k|}{2(N+1)\epsilon_k M_k}sin(\beta_s+\theta_{sl})\nonumber
\end{align}
We omit $k$ for convenience.
It is also important to observe that even if $M_k=o(\epsilon_k)$ we still have $M_k\sim \max_{s}|p_s^k-p_l^k|/\epsilon_k$. Since each $|p_l^k|=E$, an upper bound for $M_k$ is 
\begin{equation}\label{small-M}
M_k\le C\mu_k\epsilon_k+C\delta_k^2\epsilon_k^{-1}.
\end{equation}

Equation (\ref{c-12}) gives us a key observation: $|c_{1,s,l}|+|c_{2,s,l}|\sim |p_s^k-p_l^k|/(\epsilon_k M_k)$. So whenever $|c_{1,s,l}|+|c_{2,s,l}|\neq 0$ we have
$\frac{|p_s^k-p_l^k|}{\epsilon_k}\sim M_k$. In other words for each $l$, $M_{l,k}\sim \max_{t\neq l}\frac{|p_t^k-p_l^k|}{\epsilon_k}$.  Hence for any $t$, if $\frac{|p_t^k-p_l^k|}{\epsilon_k}\sim M_k$, let $M_{t,k}$ be the maximum of $|v_k-V_{t,k}|$, we have $M_{t,k}\sim M_k$. If all $\frac{|p_t^k-p_l^k|}{\epsilon_k}\sim M_k$ (\ref{M-comp}) is proved. So we prove that even if some $p_t^k$ is very close to $p_l^k$, $M_t^k$ is still comparable to $M_k$. The reason is there exists $q$ such that $\frac{|p_l^k-p_q^k|}{\epsilon_k}\sim M_k$, if $\frac{|p_t^k-p_l^k|}{\epsilon_k}=o(1)M_k$,
$$|p_t^k-p_q^k|\ge |p_l^k-p_q^k|-|p_t^k-p_l^k|\ge \frac 12 |p_l^k-p_q^k|.$$
Thus $\frac{|p_t^k-p_q^k|}{\epsilon_k}\sim M_k$ and $M_t^k\sim M_k$. (\ref{M-comp}) is established.  From now on for convenience we shall just use $M_k$. Since 
$M_k\sim \max_{s,t}|p_s^k-p_t^k|/\epsilon_k$, an upper bound of $M_k$ is
\begin{equation}\label{upper-bound-Mk}
M_k\le C\mu_k\epsilon_k. 
\end{equation}

Set $w_{l,k}=(v_k-V_{l,k})$, then we have $w_{l,k}(Q_l^k)=|\nabla w_{l,k}(Q_l^k)|=0$. Correspondingly we set
$$\tilde w_{l,k}=w_{l,k}/M_k. $$
The equation of $w_{l,k}$ can be written as 
\begin{align}\label{wlk-bs}
   &\Delta w_{l,k}+|y|^{2N}\mathfrak{h}_k(\delta_k Q_l)e^{\xi_l}w_{l,k}\\
=&-\delta_k\nabla\mathfrak{h}_k(\delta_kQ_l)(y-Q_l)|y|^{2N}e^{V_{l,k}}-\delta_k^2
\sum_{|\alpha |=2}\frac{\partial^{\alpha}\mathfrak{h}_k(\delta_kQ_l)}{\alpha !}(y-Q_l)^{\alpha}|y|^{2N}e^{V_k}\nonumber\\
&+O(\delta_k^3)|y-Q_l|^3|y|^{2N}e^{V_k} \nonumber 
\end{align}
where we omitted $k$ in $Q_l$ and $\xi_l$. $\xi_l$ comes from the Mean Value Theorem and satisfies
\begin{equation}\label{around-ls}
e^{\xi_l}=e^{V_{l,k}}(1+\frac 12w_{l,k}+O(w_{l,k}^2)). 
\end{equation}
The function $\tilde w_{l,k}$ satisfies
\begin{equation}\label{aroud-s-1}
\lim_{k\to \infty}\tilde w_{l,k}(Q_s^k+\epsilon_k z)=\frac{c_{1,s,l}z_1+c_{2,s,l}z_2}{1+\frac 18 |z|^2}
\end{equation}
and around each $Q_s^k$ (\ref{Q-bad}) holds with $M_{s,k}$ replaced by $M_k$.

Now for $|y|\sim 1$, we use $w_{l,k}(Q_l^k)=0$ to write $w_{l,k}(y)$ as
\begin{align*}
    w_{l,k}(y)&=\int_{\Omega_k}(G_k(y,\eta)-G_k(Q_l,\eta))\bigg (\mathfrak{h}_k(\delta_k Q_l)|\eta |^{2N}e^{\xi_l}w_{l,k}(\eta)\\
    &+\delta_k \nabla \mathfrak{h}_k(\delta_k Q_l)(\eta-Q_l)|\eta |^{2N}e^{V_{l,k}}\\
    &+\delta_k^2
\sum_{|\alpha |=2}\frac{\partial^{\alpha}\mathfrak{h}_k(\delta_kQ_l)}{\alpha !}(\eta-Q_l)^{\alpha}|y|^{2N}e^{V_k}\bigg )
+O(\delta_k^{N+2}). 
\end{align*}
Note that the last term is $O(\delta_k^{N+2})$  because it comes from
the oscillation of $w_{l,k}$ on $\partial \Omega_k$. The harmonic function defined by the boundary value of $w_{l,k}$ has an oscillation of  $O(\delta_k^{N+1})$ on $\partial \Omega_k$. The oscillation of this harmonic function in $B_R$ (for any fixed $R>1$) is $O(\delta_k^{N+2})$.
The regular part of the Green's function brings little error in the computation, so we have
 
 \begin{align*}
& \tilde w_{l,k}(y)\\
=&-\frac{1}{2\pi}\int_{\Omega_k}\log \frac{|y-\eta |}{|Q_l^k-\eta |}\bigg ( \tilde w_{l,k}(\eta )\mathfrak{h}_k(\delta_k Q_l^k)|\eta |^{2N}e^{\xi_l}\\
& +\sigma_k\nabla \mathfrak{h}_k(\delta_k Q_l^k)(\eta-Q_l^k)
|\eta |^{2N}e^{V_{l,k}}+\frac{\delta_k^2}{M_k}\sum_{|\alpha |=2}\frac{\partial^{\alpha}\mathfrak{h}_k(\delta_k Q_l^k)}{\alpha !}(\eta -Q_l^k)^{\alpha}|\eta |^{2N}e^{V_{l,k}}\bigg )d\eta\\
&+o(\epsilon_k)
\end{align*}
for $|y|\sim 1$.

 Around each $Q_s^k$ the $e^{\xi_l}$ can be replaced by $e^{V_{s,k}}$ with controllable error (based on Lemma \ref{t-w-1-better} and (\ref{upper-bound-Mk})). 
In order to evaluate the expression of $\tilde w_{l,k}$ we need the following identity based on (\ref{late-1})
\begin{equation}\label{late-2-r}
\int_{B(Q_s^k,\tau)}(\tilde w_{l,k}(\eta)\mathfrak{h}_k(\delta_k Q_l^k)|\eta |^{2N}e^{V_{s,k}}
+\sigma_{k}\nabla \mathfrak{h}_k(\delta_k Q_l^k)(\eta -Q_l)|\eta |^{2N}e^{V_{l,k}})d\eta=O(\epsilon_k^{\sigma}) 
\end{equation}
Note that $e^{V_{l,k}}$ in the first term was replaced by $e^{V_{s,k}}$ but in the second term above this replacement is not necessary. (\ref{late-1}) is mainly used in the evaluation of the first term. 
The proof of (\ref{late-2-r}) can be found in the appendix.

Equation (\ref{late-2-r}) also leads to a more accurate estimate of $\tilde w_{l,k}$ in regions between bubbling disks. By the Green's representation formula of $\tilde w_{l,k}$ it is easy to have, for $|y|\sim 1$,

Writing the logarithmic term above as
$$\log \frac{|y-\eta |}{|Q_l^k-\eta |}=\log \frac{|y-Q_s^k|}{|Q_l^k-Q_s^k|}+(\log \frac{|y-\eta |}{|Q_l^k-\eta |}-\log \frac{|y-Q_s^k|}{|Q_l^k-Q_s^k|}),$$
then we see that the integration related to the second term is $O(\epsilon_k)$. The integration involving the first term is $O(\epsilon_k^{\sigma})+o(\frac{\delta_k^2}{M_k})$ by (\ref{late-2-r}) and the definition of $M_k$.
Therefore
$$|\tilde w_{l,k}(y)|=o(1/\mu_k),\quad y\in B_3\setminus \cup_{s=0}^NB(Q_s^k,\tau) $$
for $\sigma\in (0,1)$. Thus this extra control of $\tilde w_{l,k}$ away from bubbling disks gives a better estimate than (\ref{Q-bad}) around $Q_l^k$: Using the same argument for Lemma \ref{t-w-1-better} we have
\begin{equation}\label{much-better-s}
|\tilde w_{l,k}(Q_l^k+\epsilon_k z)|\le o(\epsilon_k)\frac{(1+|z|)}{\log (2+|z|)},\quad |z|<\tau \epsilon_k^{-1}.
\end{equation}

From the decomposition in (\ref{late-1}) and the help from (\ref{much-better-s}) we can now estimate the integral of $\tilde w_{l,k}$ more precisely. 

\begin{align}\label{late-2}
&\int_{B(Q_s^k,\tau)}\bigg ( \tilde w_{l,k}(\eta )\mathfrak{h}_k(\delta_k Q_l^k)|\eta |^{2N}e^{\xi_l}
+\sigma_k\nabla \mathfrak{h}_k(\delta_k Q_l^k)(\eta-Q_l^k)
|\eta |^{2N}e^{V_{l,k}}\\
&+\frac{\delta_k^2}{M_k}\sum_{|\alpha |=2}
\frac{\partial^{\alpha}\mathfrak{h}_k(\delta_k Q_l^k)}{\alpha !}(\eta-Q_l^k)^{\alpha}|\eta |^{2N}e^{V_{l,k}})d\eta \nonumber\\
&=\int_{B(Q_s^k,\tau)}\tilde w_{l,k}(\eta)\mathfrak{h}_k(\delta_k Q_l^k)|\eta |^{2N}e^{V_{s,k}}d\eta+8\pi\sigma_{k}\nabla \log \mathfrak{h}_k(\delta_k Q_s^k)
(Q_s^k-Q_l^k)\nonumber\\
&\quad +2\pi\frac{\delta_k^2}{M_k}\frac{\Delta \mathfrak{h}_k(\delta_k Q_s^k)}{\mathfrak{h}_k(\delta_kQ_s^k)}|Q_s^k-Q_l^k|^2+o(\epsilon_k) \nonumber\\
&=D_{s,l}^k+o(\epsilon_k).\nonumber
\end{align}
where
$$D_{s,l}^k=\frac{\pi}{(N+1)^2}\frac{|p_s^k-p_l^k|^2}{\epsilon_k^2 M_{k}^2}M_k+2\pi\frac{\delta_k^2}{M_k}\frac{\Delta \mathfrak{h}_k(\delta_k Q_s^k)}{\mathfrak{h}_k(\delta_kQ_s^k)}|Q_s^k-Q_l^k|^2$$

Let 
$$  H_{y,l}(\eta)=\frac{1}{2\pi}\log\frac{|y-\eta|}{|Q_l^k-\eta |}. $$
Then 
\begin{align*}
    \tilde w_{l,k}(y)&=-\sum_{s\neq l}H_{y,l}(Q_s)D_{s,l}^k \\
    &-\sum_{s\neq l}\int_{B(Q_s,\tau)}\bigg (\partial_1H_{y,l}(Q_s)\eta_1+\partial_2H_{y,l}(Q_s)\eta_2)(\mathfrak{h}_k(\delta_kQ_l)|\eta |^{2N}e^{\xi_l}\tilde w_{l,k}(\eta)\\
    +&\sigma_k\nabla\mathfrak{h}_k(\delta_kQ_l)(\eta-Q_l)|\eta |^{2N}e^{V_{l,k}}+O(\delta_k^2M_k^{-1})|\eta-Q_l|^2|\eta |^{2N}e^{V_{l,k}}\bigg )d\eta+o(\epsilon_k).
\end{align*} 

After evaluation we have
\begin{align}\label{green-wlk}
    \tilde w_{l,k}(y)&=-\frac{1}{2\pi}\sum_{s\neq l}\log\frac{|y-Q_s^k|}{|Q_l^k-Q_s^k|}D_{s,l}^k
   -\sum_{s\neq l}8(\frac{y_1-Q_s^1}{|y-Q_s|^2}-\frac{Q_l^1-Q_s^1}{|Q_l-Q_s|^2})c_{1,s,l}\epsilon_k\\
   &-8(\frac{y_2-Q_s^2}{|y-Q_s|^2}-\frac{Q_l^2-Q_s^2}{|Q_s-Q_l|^2})c_{2,s,l}\epsilon_k+o(\epsilon_k). \nonumber
\end{align}
where we used
$$\int_{\mathbb R^2}\frac{z_1^2}{(1+\frac 18|z|^2)^3}dz=\int_{\mathbb R^2}\frac{z_2^2}{(1+\frac 18|z|^2)^3}dz=16\pi. $$
Recall that 
$$c_{1,s,l}=\frac{|p_s-p_l|}{2(N+1)M_k\epsilon_k}\cos (\beta_s+\theta_{sl})$$
$$c_{2,s,l}=\frac{|p_s-p_l|}{2(N+1)M_k\epsilon_k}\sin (\beta_s+\theta_{sl})$$

For $|y|\sim 1$ but away from the $N+1$ bubbling disks, we have, for $l\neq s$, 
$$v_k(y)=V_{l,k}(y)+M_k\tilde w_{l,k}(y) $$
and
$$v_k(y)=V_{s,k}(y)+M_k\tilde w_{s,k}(y). $$
Thus for $s\neq l$ we have
\begin{equation}\label{compare-10}
\frac{V_{s,k}(y)-V_{l,k}(y)}{M_k}=\tilde w_{l,k}(y)-\tilde w_{s,k}(y). \end{equation}

 In (\ref{late-1}) we consider $|z|\sim \epsilon_k^{-1}$, then we see that if 
\[|\frac{\mu_l^k-\mu_s^k}{M_k}+2\sigma_k\nabla \log \mathfrak{h}_k(\delta_kQ_s^k)(Q_l^k-Q_s^k)|\ge C\epsilon_k,\]
for a large $C$,
 it is easy to see that  (\ref{compare-10}) does not hold. So we have 
\[|\frac{\mu_l^k-\mu_s^k}{M_k}+2\sigma_k\nabla \log \mathfrak{h}_k(\delta_kQ_s^k)(Q_l^k-Q_s^k)|=O(\epsilon_k),\]
and we focus on the 
leading term $\phi_3$, which gives, for $|y|\sim 1$ away from bubbling disks, 
\[\frac{V_{s,k}(y)-V_{l,k}(y)}{M_k}
=D_{\mu}^k\epsilon_k+4Re(\frac{y^{N+1}-1}{|y^{N+1}-1|^2}\frac{\bar p_s-\bar p_l}{M_k\epsilon_k}) \epsilon_k+O(|y-Q_s^k|^2)\epsilon_k+o(\epsilon_k), \]
where 
\[D_{\mu}^k=\frac{\mu_l^k-\mu_s^k}{M_k\epsilon_k}+2\frac{\sigma_k}{\epsilon_k}\nabla \log \mathfrak{h}_k(\delta_kQ_s^k)(Q_l^k-Q_s^k).\]
On the other hand, for $y\in B_5\setminus (\cup_{t=1}^NB(Q_t^k,\tau_1))$, 
\begin{align*}
 &\tilde w_{l,k}(y)-\tilde w_{s,k}(y)\\
 =&-\frac{1}{2\pi}\sum_{m,m\neq l}\log \frac{|y-Q_m|}{|Q_l-Q_m|}D_{m,l}^k\\
 &+8\epsilon_k\sum_{m,m\neq l}\bigg ((\frac{y_1-Q_m^1}{|y-Q_m|^2}-\frac{Q_l^1-Q_m^1}{|Q_l-Q_m|^2})\frac{|p_m-p_l|}{2(N+1)M_k\epsilon_k}\cos(\beta_m+\theta_{ml})\\
 &+(\frac{y_2-Q_m^2}{|y-Q_m|^2}-\frac{Q_l^2-Q_m^2}{|Q_l-Q_m|^2})\frac{|p_m-p_l|}{2(N+1)M_k\epsilon_k}\sin (\beta_m+\theta_{ml})\bigg )\\
 &+\frac{1}{2\pi}\sum_{m,m\neq s}\log \frac{|y-Q_m|}{|Q_s-Q_m|}D_{m,s}^k\\
 &-8\epsilon_k\sum_{m,m\neq s}\bigg ((\frac{y_1-Q_m^1}{|y-Q_m|^2}-\frac{Q_s^1-Q_m^1}{|Q_s-Q_m|^2})\frac{|p_m-p_s|}{2(N+1)M_k\epsilon_k}\cos(\beta_m+\theta_{ms})\\
 &+(\frac{y_2-Q_m^2}{|y-Q_m|^2}-\frac{Q_s^2-Q_m^2}{|Q_s-Q_m|^2})\frac{|p_m-p_s|}{2(N+1)M_k\epsilon_k}\sin (\beta_m+\theta_{ms})\bigg )
\end{align*}
for all $l\neq s$. If we fix a set of $l,s$ that corresponds to the largest $|D_{s,l}^k|$ and we consider $y$ close to $Q_s^k$. If we use $y=e^{i\beta_s}+z$ by abusing the notation $z$, then we have
\[y^{N+1}=(e^{i\beta_s}(1+ze^{-i\beta_s}))^{N+1}=1+(N+1)ze^{-i\beta_s}+O(|z|^2).\]
Therefore
\begin{align*}
&4Re(\frac{y^{N+1}-1}{|y^{N+1}-1|^2}\frac{\bar p_s-\bar p_l}{M_k\epsilon_k}\\
=&\frac{4|p_s-p_l|}{(N+1)|z|^2M_k\epsilon_k}\bigg (z_1\cos(\beta_s+\beta_{sl})+z_2\sin(\beta_s+\beta_{sl})+O(|z|^2)\bigg ).
\end{align*}
In the expression of $\tilde w_{l,k}(y)-\tilde w_{s,k}(y)$, we identify the leading term, which is 
\begin{align*}
8\epsilon_k\bigg ((\frac{y_1-Q_s^1}{|y-Q_s|^2}-\frac{Q_l^1-Q_s^1}{|Q_l-Q_s|^2})\frac{|p_s-p_l|}{2(N+1)M_k\epsilon_k}\cos(\beta_s+\theta_{sl})\\
+(\frac{y_2-Q_s^2}{|y-Q_s|^2}-\frac{Q_l^2-Q_s^2}{|Q_l-Q_s|^2})\frac{|p_s-p_l|}{2(N+1)M_k\epsilon_k}\sin(\beta_s+\theta_{sl})\bigg )\\
-\frac{1}{2\pi}\log \frac{|y-Q_s|}{|Q_l-Q_s|}D_{s,l}^k.
\end{align*}
If we use $y=e^{i\beta_s}+z$ for $|z|$ small and replace $Q_s$ by $e^{i\beta_s}$ because their difference is $o(\epsilon_k)$. Then the expression above has this leading term:
\[\frac{4|p_s-p_l|}{(N+1)|z|^2M_k\epsilon_k}\bigg (z_1\cos(\beta_s+\beta_{sl})+z_2\sin(\beta_s+\beta_{sl})\bigg )-\frac 1{2\pi}\log \frac{|z|}{| e^{i\beta_l}-e^{i\beta_s}|}D_{s,l}^k.\]
Thus we obtain $D_{s,l}^k/\epsilon_k=o(1)$. Therefore 
\begin{equation}\label{dsl0}
D_{s,l}^k=o(\epsilon_k),\quad \forall s\neq l.
\end{equation}
With this updated information we write $\tilde w_{l,k}(y)-\tilde w_{s,k}(y)$ as
\begin{align*}
 &\tilde w_{l,k}(y)-\tilde w_{s,k}(y)\\
 &=8\epsilon_k\sum_{m,m\neq l}\bigg ((\frac{y_1-\cos \beta_m}{|y-e^{i\beta_m}|^2}-\frac{\cos \beta_l-\cos \beta_m}{|e^{i\beta_l}-e^{i\beta_m}|^2})\frac{|p_m-p_l|}{2(N+1)M_k\epsilon_k}\cos(\beta_m+\theta_{ml})\\
 &+(\frac{y_2-\sin \beta_m}{|y-e^{i\beta_m}|^2}-\frac{\sin\beta_l-\sin \beta_m}{|e^{i\beta_l}-e^{i\beta_m}|^2})\frac{|p_m-p_l|}{2(N+1)M_k\epsilon_k}\sin (\beta_m+\theta_{ml})\bigg )\\
 &-8\epsilon_k\sum_{m,m\neq s}\bigg ((\frac{y_1-\cos \beta_m}{|y-e^{i\beta_m}|^2}-\frac{\cos \beta_s-\cos \beta_m}{|e^{i\beta_s}-e^{i\beta_m}|^2})\frac{|p_m-p_s|}{2(N+1)M_k\epsilon_k}\cos(\beta_m+\theta_{ms})\\
 &+(\frac{y_2-\sin \beta_m}{|y-e^{i\beta_m}|^2}-\frac{\sin \beta_s-\sin \beta_m}{|e^{i\beta_s}-e^{i\beta_m}|^2})\frac{|p_m-p_s|}{2(N+1)M_k\epsilon_k}\sin (\beta_m+\theta_{ms})\bigg )+o(\epsilon_k)
\end{align*}
for $|y|\in B_5\setminus (\cup_{l=1}^N B(Q_l^k,\tau_1))$. 

\medskip

Now in particular we take $l=0$ and we use the following notations: $\tilde w_k$, $V_k$, $c_{1,s}$, $c_{2,s}$, $\theta_s$, instead of $\tilde w_{0}^k$, $v_{0,k}$, $c_{1,s,0}$, $c_{2,s,0}$, $\theta_{s,0}$. 

The expression of $\tilde w_k$ (see (\ref{green-wlk}) for example) gives
\begin{align*}
&\nabla \tilde w_k(y)\\
=&\int_{\Omega_k} \nabla_y G(y,\eta)\bigg (\mathfrak{h}_k(\delta_ke_1)|\eta |^{2N}e^{\xi_k}\tilde w_k(\eta)
+\sigma_k\nabla\mathfrak{h}_k(\delta_ke_1)(\eta-e_1)|\eta |^{2N}e^{V_k(\eta)}\\
&+\frac{\delta_k^2}{M_k}\sum_{|\alpha |=2}\frac{\partial^{\alpha}\mathfrak{h}_k(\delta_ke_1)(\eta -e_1)^{\alpha}}{\alpha !}|\eta |^2{2N}e^{V_k(\eta)}\bigg )d\eta
+o(\epsilon_k), 
\end{align*}
for $y\in B_5\setminus (\cup_{s=1}^N B(Q_s^k,\tau_1))$.
Now we take $y=e_1$, we have 
\begin{align*}
&0=\nabla \tilde w_k(e_1)\\
=&\int_{\Omega_k} (-\frac{1}{2\pi})\frac{e_1-\eta}{|e_1-\eta |^2}\bigg (\mathfrak{h}_k(\delta_ke_1)|\eta |^{2N}e^{\xi_k}\tilde w_k(\eta)
+\sigma_k\nabla\mathfrak{h}_k(\delta_ke_1)(\eta-e_1)|\eta |^{2N}e^{V_k(\eta)}\\
&+\frac{\delta_k^2}{M_k}\sum_{|\alpha |=2}\frac{\partial^{\alpha}\mathfrak{h}_k(\delta_ke_1)(\eta -e_1)^{\alpha}}{\alpha !}|\eta |^2{2N}e^{V_k(\eta)}\bigg )d\eta
+o(\epsilon_k), 
\end{align*}
Obviously we will integral each of the two components in $B(Q_s^k,\tau_1)$ for $\tau_1>0$ small. Then we observe from (\ref{dsl0}) 
that
\begin{align*}\int_{B(Q_l,\tau_1)}\bigg (\mathfrak{h}_k(\delta_ke_1)|\eta |^{2N}e^{\xi_k}\tilde w_k(\eta)
+\sigma_k\nabla\mathfrak{h}_k(\delta_ke_1)(\eta-e_1)|\eta |^{2N}e^{V_k(\eta)}\\
+\frac{\delta_k^2}{M_k}\sum_{|\alpha |=2}\frac{\partial^{\alpha}\mathfrak{h}_k(\delta_ke_1)(\eta -e_1)^{\alpha}}{\alpha !}|\eta |^2{2N}e^{V_k(\eta)}\bigg )d\eta=o(\epsilon_k).\end{align*}
Based on this we use the following format: If $f$ is a smooth function,
\begin{align*}
&\int_{B(Q_s,\tau_1)} f(\eta) \bigg (\mathfrak{h}_k(\delta_ke_1)|\eta |^{2N}e^{\xi_k}\tilde w_k(\eta)
+\sigma_k\nabla\mathfrak{h}_k(\delta_ke_1)(\eta-e_1)|\eta |^{2N}e^{V_k(\eta)}\\
&\qquad +\frac{\delta_k^2}{M_k}\sum_{|\alpha |=2}\frac{\partial^{\alpha}\mathfrak{h}_k(\delta_ke_1)(\eta -e_1)^{\alpha}}{\alpha !}|\eta |^2{2N}e^{V_k(\eta)}\bigg )d\eta\\
&=\partial_1f(e^{i\beta_s})c_{1s}\cdot 16\pi \epsilon_k+\partial_1 f(e^{i\beta_s})c_{2s}\cdot 16\pi \epsilon_k +o(\epsilon_k). 
\end{align*}
Then we replace $f(\eta_1,\eta_2)$ by \[f_1(\eta_1,\eta_2)=(-\frac{1}{2\pi})\frac{1-\eta_1}{(1-\eta_1)^2+\eta_2^2}\] and 
\[f_2(\eta_1,\eta_2)=\frac{1}{2\pi}\frac{\eta_2}{(1-\eta_1)^2+\eta_2^2}.\]
Then we have, from the expressions of $c_{1,s}$, $c_{2,s}$ in (\ref{c-12}), that 
\begin{align*}
0&=\partial_1\tilde w_k(e_1)\\
&=16\pi\epsilon_k\sum_{s=1}^N\bigg (\frac{\cos(\beta_s+\theta_s)\cos \beta_s+\sin \beta_s\sin (\beta_s+\theta_s)}{4\pi(1-\cos \beta_s)}\frac{|p_s|}{2(N+1)M_k\epsilon_k}\bigg )+o(\epsilon_k)\\
&=4\epsilon_k\sum_{s=1}^N\frac{\cos \theta_s}{1-\cos \beta_s}\frac{|p_s|}{2(N+1)M_k\epsilon_k} +o(\epsilon_k). 
\end{align*}
Similarly
\begin{align*}
0&=\partial_2\tilde w_k(e_1)\\
&=16\pi\epsilon_k\sum_{s=1}^N\bigg (\frac{\cos(\beta_s+\theta_s)\sin \beta_s-\cos \beta_s\sin (\beta_s+\theta_s)}{4\pi(1-\cos \beta_s)}\frac{|p_s|}{2(N+1)M_k\epsilon_k}\bigg )+o(\epsilon_k)\\
&=-4\epsilon_k\sum_{s=1}^N\frac{\sin \theta_s}{1-\cos \beta_s}\frac{|p_s|}{2(N+1)M_k\epsilon_k} +o(\epsilon_k). 
\end{align*}
If we use $a_s$ to denote
\[a_s:=\lim_{k\to \infty}\frac{|p_s^k|}{(1-cos\beta_s)M_k\epsilon_k},\quad s=1,...,N,\]
then we have
\begin{equation}\label{eq-ar-1}
\sum_{s=1}^Na_s\cos \theta_s=0
\end{equation}
\begin{equation}\label{eq-ar-2}
\sum_{s=1}^Na_s\sin \theta_s=0,
\end{equation}
where $a_s\ge 0$ for all $s$. Taking the sum of the squares of (\ref{eq-ar-1}) and (\ref{eq-ar-2}) we obtain 
\[\sum_{s=1}^Na_s^2+\sum_{s<t}2a_sa_t\cos(\theta_s-\theta_t)=0. \]
Since $\sum_sa_s^2\ge 2\sum_{s<t}a_sa_t$, we have
\[\sum_{s<t}2a_sa_t(1+\cos(\theta_s-\theta_t))\le 0.\]
Since each term on the left is obviously non-negative, we know each 
\[a_sa_t(1+\cos (\theta_s-\theta_t))=0,\quad \forall s<t. \] 
If there is only one $a_s>0$, it is easy to see that (\ref{eq-ar-1}) and (\ref{eq-ar-2}) cannot both hold. If there are three $a_t's>0$, it is also elementary to see this is not possible: say $a_1,a_2,a_3>0$, then they have to be equal. Then we see that we must have 
\[\theta_1-\theta_2=\pm \pi,\quad \theta_2-\theta_3=\pm \pi,\quad \theta_1-\theta_3=\pm \pi .\] 
Obviously these three equations cannot hold at the same time. So the only situation left is there are exactly two $a_t's$ positive. All other $a_ts$ are zero. Since $p_0=0$, this means there are exactly two $p_{s_1}^k$, $p_{s_2}^k$ such that 
\begin{equation}\label{p-ex-1}\lim_{k\to \infty}\frac{p_{s_1}^k}{\epsilon_kM_k}=-\lim_{k\to \infty}\frac{p_{s_2}^k}{\epsilon_k M_k}\neq 0,\quad \lim_{k\to \infty}\frac{p_t^k}{\epsilon_kM_k}=0,\quad \forall t\neq s_1,s_2.
\end{equation}

If we apply the same argument to $\tilde w_l^k$. Then from $\nabla \tilde w_l^k(Q_l^k)=0$ we would get exactly $p_{l_1}^k$ and $p_{l_2}^k$ different from $p_l^k$ and 
\[\lim_{k\to \infty} \frac{p_{l_1}^k-p_l^k}{\epsilon_k M_k}=-\lim_{k\to \infty}\frac{p_{l_2}^k-p_l^k}{\epsilon_k M_k}\neq 0,\quad \lim_{k\to \infty}\frac{p_t^k-p_l^k}{\epsilon_k M_k}=0, \forall t\neq l_1,l_2.\]
Then it is easy to see that this is only possible when we have $N=2$ because if $N=1$, we would have just one $a_s\neq 0$, which is not possible based on (\ref{eq-ar-1}) and (\ref{eq-ar-2}). If $N\ge 3$, we have to have $p_t^k$ that satisfies 
\[\lim_{k\to \infty}\frac{|p_t^k|}{\epsilon_k M_k}=0,\quad \mbox{and}\quad \lim_{k\to \infty} \frac{p_t^k-p_{s_1}^k}{\epsilon_kM_k}=0\]
which is a contradiction to (\ref{p-ex-1}). 

Finally we rule out the case $N=2$. In this case we have $p_0^k=0$, 
\begin{equation}\label{p-ex-2}\lim_{k\to \infty}\frac{p_1^k}{\epsilon_kM_k}=-\lim_{k\to \infty}\frac{p_2^k}{\epsilon_kM_k}\neq 0.
\end{equation}
However from $\tilde w_1^k(p_1^k)=0$ we have 
\[\lim_{k\to \infty}\frac{p_2^k-p_1^k}{\epsilon_kM_k}=-\lim_{k\to \infty}\frac{0-p_1^k}{\epsilon_kM_k}\neq 0,\]
which is a contradiction to (\ref{p-ex-2}).
Lemma \ref{small-other} is established. $\Box$

\medskip

Proposition \ref{key-w8-8} is an immediate consequence of Lemma \ref{small-other}.  $\Box$.

\medskip

Now we finish the proof of Theorem \ref{vanish-first-h}.

Let $\hat w_k=w_k/\tilde \delta_k$. (Recall that $\tilde \delta_k=\delta_k |\nabla \mathfrak{h}_k(0)|+\delta_k^2\mu_k$). If $|\nabla \mathfrak{h}_k(0)|/(\delta_k\mu_k)\to \infty$, we see that in this case 
$\tilde \delta_k\sim \delta_k\mu_k|\nabla \mathfrak{h}_k(0)|$. The equation of $\hat w_k$ is
\begin{equation}\label{hat-w}
\Delta \hat w_k+|y|^{2N}e^{\xi_k}\hat w_k=a_k\cdot (e_1-y)|y|^{2N}e^{V_k}+b_ke^{V_k}|y-e_1|^{2}|y|^{2N},
\end{equation}
in $\Omega_k$, where $a_k=\delta_k\nabla \mathfrak{h}_k(0)/\tilde \delta_k$, $b_k=o(1/\mu_k)$. 
By Proposition \ref{key-w8-8}, $|\hat w_k(y)|\le C$. Before we carry out the remaining part of the proof we observe that $\hat w_k$ converges to a harmonic function in $\mathbb R^2$ minus finite singular points. Since $\hat w_k$ is bounded, all these singularities are removable. Thus $\hat w_k$ converges to a constant. Based on the information around $e_1$, we shall prove that this constant is $0$. However, looking at the right hand side the equation,
$$ (e_1-y)|y|^{2N}e^{V_k}\rightharpoonup \sum_{l=1}^N8\pi  (e_1-e^{i\beta_l})\delta_{e^{i\beta_l}}. $$
we will get a contradiction by comparing the Pohozaev identities of $v_k$ and $V_k$, respectively.

Now we use the notation $W_k$ again and use Proposition \ref{key-w8-8} to rewrite the equation for $W_k$.
Let
$$W_k(z)=\hat w_k(e_1+\epsilon_k z), \quad |z|< \delta_0 \epsilon_k^{-1} $$
for $\delta_0>0$ small. Then from Proposition \ref{key-w8-8} we have
\begin{equation}\label{h-exp}
\mathfrak{h}_k(\delta_ky)=\mathfrak{h}_k(\delta_k e_1)+\delta_k \nabla \mathfrak{h}_k(\delta_k e_1)(y-e_1)+O(\delta_k^2)|y-e_1|^2,
\end{equation}
\begin{equation}\label{y-1}
|y|^{2N}=|e_1+\epsilon_k z|^{2N}=1+O(\epsilon_k)|z|,
\end{equation}
\begin{equation}\label{v-radial}
V_k(e_1+\epsilon_k z)+2\log \epsilon_k=U_k(z)+O(\epsilon_k)|z|+O(\epsilon_k^2)(\log (1+|z|))^2
\end{equation}
and
\begin{equation}\label{xi-radial}
\xi_k(e_1+\epsilon_k z)+2\log \epsilon_k=U_k(z)+O(\epsilon_k)(1+|z|).
\end{equation}
Using (\ref{h-exp}),(\ref{y-1}),(\ref{v-radial}) and (\ref{xi-radial}) in (\ref{hat-w}) we write the equation of $W_k$ as
\begin{equation}\label{W-rough}
\Delta W_k+\mathfrak{h}_k(\delta_k e_1)e^{U_k(z)}W_k=-\epsilon_k a_k\cdot ze^{U_k(z)}+E_w, \quad 0<|z|<\delta_0 \epsilon_k^{-1}
\end{equation}
where
\begin{equation}\label{ew-rough-2}
E_w(z)=O(\epsilon_k)(1+|z|)^{-3}, \quad |z|<\delta_0 \epsilon_k^{-1}.
\end{equation}

Since $\hat w_k$ obviously converges to a global harmonic function with removable singularity, we have $\hat w_k\to \bar c$ for some $\bar c\in \mathbb R$. Then we claim that
\begin{lem}\label{bar-c-0} $\bar c=0$.
\end{lem}

\noindent{\bf Proof of Lemma \ref{bar-c-0}:}

 If $\bar c\neq 0$, we use $W_k(z)=\bar c+o(1)$ on $B(0,\delta_0 \epsilon_k^{-1})\setminus B(0, \frac 12\delta_0 \epsilon_k^{-1})$ and consider the projection of $W_k$ on $1$:
 $$g_0(r)=\frac 1{2\pi}\int_{0}^{2\pi}W_k(re^{i\theta})d\theta. $$
If we use $F_0$ to denote the projection to $1$ of the right hand side we have, using the rough estimate of $E_w$ in (\ref{ew-rough-2})
$$g_0''(r)+\frac 1r g_0'(r)+\mathfrak{h}_k(\delta_ke_1)e^{U_k(r)}g_0(r)=F_0,\quad 0<r<\delta_0 \epsilon_k^{-1} $$
where
$$F_0(r)=O(\epsilon_k)(1+|z|)^{-3}. $$
In addition we also have
$$\lim_{k\to \infty} g_0(\delta_0 \epsilon_k^{-1})=\bar c+o(1). $$
For simplicity we omit $k$ in some notations. By the same argument as in Lemma \ref{w-around-e1},  we have
$$g_0(r)=O(\epsilon_k)\log (2+r),\quad 0<r<\delta_0 \epsilon_k^{-1}. $$
Thus $\bar c=0$.
Lemma \ref{bar-c-0} is established. $\Box$

\medskip

Based on Lemma \ref{bar-c-0} and standard Harnack inequality for elliptic equations we have
\begin{equation}\label{small-til-w}
\tilde w_k(x)=o(1),\,\,\nabla \tilde w_k(x)=o(1),\,\, x\in B_3\setminus (\cup_{l=1}^N(B(e^{i\beta_l},\delta_0)\setminus B(e^{i\beta_l}, \delta_0/8))).
\end{equation}
Equation (\ref{small-til-w}) is equivalent to $w_k=o(\tilde \delta_k)$ and $\nabla w_k=o(\tilde \delta_k)$ in the same region.

\medskip

In the next step we consider the difference between two Pohozaev identities.
 For $s=1,...,N$ we consider the Pohozaev identity around $Q_s^k$. Let $\Omega_{s,k}=B(Q_s^k,r)$ for small $r>0$. For $v_k$ we have
\begin{align}\label{pi-vk}
\int_{\Omega_{s,k}}\partial_{\xi}(|y|^{2N}\mathfrak{h}_k(\delta_ky))e^{v_k}-\int_{\partial \Omega_{s,k}}e^{v_k}|y|^{2N}\mathfrak{h}_k(\delta_ky)(\xi\cdot \nu)\\
=\int_{\partial \Omega_{s,k}}(\partial_{\nu}v_k\partial_{\xi}v_k-\frac 12|\nabla v_k|^2(\xi\cdot \nu))dS. \nonumber
\end{align}
where $\xi$ is an arbitrary unit vector. Correspondingly the Pohozaev identity for $V_k$ is

\begin{align}\label{pi-Vk}
\int_{\Omega_{s,k}}\partial_{\xi}(|y|^{2N}\mathfrak{h}_k(\delta_ke_1))e^{V_k}-\int_{\partial \Omega_{s,k}}e^{V_k}|y|^{2N}\mathfrak{h}_k(\delta_k e_1)(\xi\cdot \nu)\\
=\int_{\partial \Omega_{s,k}}(\partial_{\nu}V_k\partial_{\xi}V_k-\frac 12|\nabla V_k|^2(\xi\cdot \nu))dS. \nonumber
\end{align}
Using $w_k=v_k-V_k$ and $|w_k(y)|\le C\tilde \delta_k$ we have
\begin{align*}
&\int_{\partial \Omega_{s,k}}(\partial_{\nu}v_k\partial_{\xi}v_k-\frac 12|\nabla v_k|^2(\xi\cdot \nu))dS\\
=&\int_{\partial \Omega_{s,k}}(\partial_{\nu}V_k\partial_{\xi}V_k-\frac 12|\nabla V_k|^2(\xi\cdot \nu))dS\\
&+\int_{\partial \Omega_{s,k}}(\partial_{\nu}V_k\partial_{\xi}w_k+\partial_{\nu}w_k\partial_{\xi}V_k-(\nabla V_k\cdot \nabla w_k)(\xi\cdot \nu))dS+o(\tilde \delta_k).
\end{align*}
If we just use crude estimate: $\nabla w_k=o(\tilde \delta_k)$, we have
\begin{align*}&\int_{\partial \Omega_{s,k}}(\partial_{\nu}v_k\partial_{\xi}v_k-\frac 12|\nabla v_k|^2(\xi\cdot \nu))dS\\
-&\int_{\partial \Omega_{s,k}}(\partial_{\nu}V_k\partial_{\xi}V_k-\frac 12|\nabla V_k|^2(\xi\cdot \nu))dS
=o(\tilde \delta_k).
\end{align*}
The difference on the second terms is minor: If we use the expansion of $v_k=V_k+w_k$ and that of $\mathfrak{h}_k(\delta_ky)$ around $e_1$, it is easy to obtain
$$\int_{\partial \Omega_{s,k}}e^{v_k}|y|^{2N}\mathfrak{h}_k(\delta_ky)(\xi\cdot \nu)-\int_{\partial \Omega_{s,k}}e^{V_k}|y|^{2N}\mathfrak{h}_k(\delta_ke_1)(\xi\cdot \nu)=o(\tilde \delta_k). $$
To evaluate the first term,  we use
\begin{align}\label{imp-1}
&\partial_{\xi}(|y|^{2N}\mathfrak{h}_k(\delta_ky))e^{v_k}\\
=&\partial_{\xi}(|y|^{2N}\mathfrak{h}_k(\delta_ke_1)+|y|^{2N}\delta_k\nabla \mathfrak{h}_k(\delta_ke_1)(y-e_1)+O(\delta_k^2))e^{V_k}(1+w_k+O(\delta_k^2\mu_k))\nonumber\\
=&\partial_{\xi}(|y|^{2N})\mathfrak{h}_k(\delta_k e_1)e^{V_k}+\delta_k\partial_{\xi}(|y|^{2N}\nabla \mathfrak{h}_k(\delta_ke_1)(y-e_1))e^{V_k}\nonumber\\
&+\partial_{\xi}(|y|^{2N}\mathfrak{h}_k(\delta_ke_1))e^{V_k}w_k+O(\delta_k^2\mu_k)e^{V_k}.\nonumber
\end{align}

For the third term on the right hand side of (\ref{imp-1}) we use the equation for $w_k$:
$$\Delta w_k+\mathfrak{h}_k(\delta_ke_1)e^{V_k}|y|^{2N}w_k=-\delta_k\nabla \mathfrak{h}_k(\delta_ke_1)\cdot (y-e_1)|y|^{2N}e^{V_k}+O(\delta_k^2)e^{V_k}|y|^{2N}.$$
From integration by parts we have
\begin{align}\label{extra-1}
&\int_{\Omega_{s,k}}\partial_{\xi}(|y|^{2N})\mathfrak{h}_k(\delta_ke_1)e^{V_k}w_k\nonumber\\
=&2N\int_{\Omega_{s,k}}|y|^{2N-2}y_{\xi}\mathfrak{h}_k(\delta_ke_1)e^{V_k}w_k\nonumber\\
=&2N\int_{\Omega_{s,k}}\frac{y_{\xi}}{|y|^2}(-\Delta w_k-\delta_k \nabla\mathfrak{h}_k(\delta_ke_1)(y-e_1)|y|^{2N}e^{V_k}+O(\delta_k^2)e^{V_k}|y|^{2N})\nonumber\\
=&-2N\delta_k\int_{\Omega_{s,k}}\frac{y_{\xi}}{|y|^{2}}\nabla \mathfrak{h}_k(\delta_ke_1)(y-e_1)|y|^{2N}e^{V_k}\nonumber\\
&+2N\int_{\partial \Omega_{s,k}}(\partial_{\nu}(\frac{y_{\xi}}{|y|^2})w_k-\partial_{\nu}w_k\frac{y_{\xi}}{|y|^2})+o(\tilde \delta_k)\nonumber\\
=&\nabla\mathfrak{h}_k(\delta_ke_1)\bigg (-16N\delta_k\pi(e^{i\beta_s}\cdot \xi)(e^{i\beta_s}-e_1)+O(\mu_k\epsilon_k^2)\bigg )+o(\tilde \delta_k),
\end{align}
where we have used $\nabla w_k, w_k=o(\tilde \delta_k)$ on $\partial \Omega_{s,k}$.
For the second term on the right hand side of (\ref{imp-1}), we have
\begin{align}\label{imp-2}
&\int_{\Omega_{s,k}}\delta_k\partial_{\xi}(|y|^{2N}\nabla \mathfrak{h}_k(\delta_ke_1)(y-e_1))e^{V_k}\\
=&2N\delta_k\int_{\Omega_{s,k}}y_{\xi}|y|^{2N-2}\nabla \mathfrak{h}_k(\delta_ke_1)(y-e_1)e^{V_k}+
\delta_k\int_{\Omega_{s,k}}|y|^{2N}\partial_{\xi}\mathfrak{h}_k(\delta_ke_1)e^{V_k} \nonumber \\
=&\nabla \mathfrak{h}_k(\delta_ke_1)\big (16N\pi\delta_k(e^{i\beta_s}\cdot \xi)(e^{i\beta_s}-e_1)+O(\mu_k\epsilon_k^2)\big )\nonumber\\
&+
\delta_k\partial_{\xi}\mathfrak{h}_k(\delta_ke_1)(8\pi+O(\mu_k\epsilon_k^2))+o(\tilde \delta_k).\nonumber
\end{align}
Using (\ref{extra-1}) and (\ref{imp-2}) in the difference between (\ref{pi-vk}) and (\ref{pi-Vk}), we have
$$\delta_k \partial_{\xi}\mathfrak{h}_k(\delta_ke_1)(1+O(\mu_k\epsilon_k^2))=o(\tilde \delta_k). $$
Thus $\nabla \mathfrak{h}_k(\delta_ke_1)=O(\delta_k\mu_k)$. Theorem \ref{vanish-first-h} is established.  $\Box$

\section{Proof of Theorem \ref{main-thm-1}}

First we prove the case $N\ge 2$. In \cite{wei-zhang-adv} we have already proved 
that 
$$\Delta (\log \mathfrak{h}_k)(0)=O(\delta_k^{-2}\mu_ke^{-\mu_k})+O(\delta_k). $$
So if $\delta_k/( \mu_k^{\frac 12}\epsilon_k)\to \infty$ there is nothing to prove. So we only consider the case that $\delta_k\le C\mu_k^{\frac 12}\epsilon_k$. By in this case $\epsilon_k^{-1}\delta_k^2\le C\epsilon_k^{\epsilon}$ for some $\epsilon\in (0,1)$. The whole argument of Proposition \ref{key-w8-8} can be
employed to prove
\begin{equation}\label{2nd-w}
|w_k(y)|\le C\delta_k^2\mu_k^{\frac{7}{4}} 
\end{equation}
In order to employ the same strategy of proof, one needs to have three things: first $\epsilon_k^{-1}\delta_k^2=O(\epsilon_k^{\epsilon})$. This is clear from the definition of $\delta_k$. Second, in the proof of Lemma \ref{small-other} we need 
$$O(\delta_k^{N+2}/M_k)=o(\epsilon_k), $$
where $M_k>\delta_k^2\mu_k^{7/4}$. Since $\delta_k\le C\mu_k^{\frac 12}\epsilon_k$ and $N\ge 1$, the required inequality holds. Thirdly, we need to have
$$\frac{\delta_k^3}{M_k}=o(\epsilon_k).$$ This is used in (\ref{subtle-1}). From the requirement on $\delta_k$ and the definition of $M_k$  This clearly also holds. The proof of Proposition \ref{key-w8-8} follows. Thus for $N\ge 2$ we also have (\ref{2nd-w}). 

\medskip
 
The precise upper bound of $w_k$ in (\ref{2nd-w}) leads to the vanishing rate of the Laplacian estimate for $N\ge 2$ and some cases of $N=1$: If we use 
$$W_k(z)=w_k(e_l+\epsilon_kz)/(\delta_k^2\mu_k^{\frac 74}), \quad |z|<\tau \epsilon_k^{-1}$$ where $e_l\neq e_1$. 
We shall show that the projection of $W_k$ over $1$ is not bounded when $|z|\sim \epsilon_k^{-1}$, which gives the desired contradiction. 

We write the equation of $w_k$ as
$$\Delta w_k+|y|^{2N}e^{\xi_k}w_k
=(\mathfrak{h}_k(\delta_k e_1)-\mathfrak{h}_k(\delta_k y))|y|^{2N}e^{V_k}. 
$$
Then for $l\neq 1$, 
\begin{align*}
&\Delta W_k(z)+e^{U_k}W_k(z)\\
=&a_0e^{U_k}+a_1ze^{U_k}+\frac 1{2\mu_k^{7/4}}\Delta \mathfrak{h}_k(0)|z|^2e^{U_k}+\frac{1}{\mu_k^{7/4}}R_2(\theta)|z|^2e^{U_k}
+O(\epsilon_k^{\epsilon}(1+|z|)^{-3}).
\end{align*}
where 
$$a_0=(\mathfrak{h}_k(\delta_ke_1)-\mathfrak{h}_k(\delta_ke_l))/(\delta_k^2\mu_k^{7/4}),$$
$$a_1=-\nabla\mathfrak{h}_k(\delta_ke_l)/(\delta_k\mu_k^{7/4}),$$
$R_2$ is the collection of spherical harmonic functions of degree $2$. Note that the assumption $l\neq 1$ means there is no appearance of $\epsilon_k$ or $\epsilon_k^2$ in the equation for $W_k$.

Let $g_k(r)$ be the projection of $W_k$ on $1$, by the same ODE analysis as before, we see that $g_k$ satisfies
$$g_k''+\frac 1rg_k'(r)+e^{U_k}g_k=E_k $$
where
$$E_k(r)=O(\epsilon_k^{\epsilon})(1+r)^{-3}+\frac 1{2\mu_k^{7/4}}\Delta (\log \mathfrak{h}_k)(0)r^2e^{U_k}. $$

Using the same argument as in Lemma \ref{w-around-e1}, we have
 $$g_k(r)\sim \Delta (\log \mathfrak{h}_k)(0)
(\log r)^2 \mu_k^{-7/4},\quad r>10. $$

Clearly if $\Delta (\log \mathfrak{h}_k(0)\neq 0$ we obtain a violation of the bound of $w_k$ for $r\sim \epsilon_k^{-1}$. Theorem \ref{main-thm-1} for $N\ge 2$ is proved under the assumption 
\begin{equation}\label{assump-small}
\epsilon_k^{-1}|Q_s^k-e^{i\beta_s}|\le \epsilon_k^{\epsilon},\quad s=1,...,N.
\end{equation}
We need this assumption because the $\xi_k$ function that comes from the equation of $w_k$ needs to tend to $U$ after scaling. From (3.13) in \cite{wei-zhang-adv}, 
$|Q_s^k-e^{i\beta_s}|=O(\delta_k^2)+O(\mu_ke^{-\mu_k})$.
If $\delta_k^2\epsilon_k^{-1}\ge C$, the argument in Theorem \ref{vanish-first-h} cannot be used because either $\xi_k$ does not tend $U$ or $c_0=0$ cannot be proved.
For $N\ge 2$, this is not a problem because we only consider $\delta_k\le C\mu_k^{\frac 12}\epsilon_k$. 

\medskip

Next we prove Theorem \ref{main-thm-1} for $N=1$ and $\delta_k\le \mu_k\epsilon_k$. The reader could see immediately that the same proof for the case $N\ge 2$ still works. 
So the only remaining case is

\medskip

\noindent{\bf Proof of Theorem \ref{main-thm-1} for $N=1$ and $\delta_k\ge \mu_k\epsilon_k$. }

In this case we  write the equation of $w_k$ as
$$\Delta w_k+|y|^2\mathfrak{h}_k(\delta_ky)e^{v_k}-|y|^2\mathfrak{h}_k(\delta_ke_1)e^{V_k}=0. $$

From $0=\nabla w_k(e_1)$ we have
\begin{equation}\label{location-mid-1}
0=\int_{\Omega_k}\nabla_1G_k(e_1,\eta)|\eta |^2(\mathfrak{h}_k(\delta_k \eta)e^{v_k}-\mathfrak{h}_k(\delta_ke_1)e^{V_k})d\eta+O(\delta_k^3)
\end{equation}
Note that $v_k$ is close to another global solution $\bar V_k$ which matches with a local maximum of $v_k$ at $Q_2^k$. 
Evaluating the right hand side of (\ref{location-mid-1}) we have
$$\nabla_1 G_k(e_1,Q_2^k)-\nabla_1 G_k(e_1,e^{i\pi})=O(\epsilon_k^2\mu_k)+O(\delta_k^3). $$
This expression gives
$$Q_2^k-e^{i\pi}=O(\delta_k^3)+O(\mu_k\epsilon_k^2). $$
This estimate will lead to a better estimate of $w_k$ outside the two bubbling disks. 
From the Green's representation for $w_k$ we now obtain 
$$w_k(y)=\int_{\Omega_k}(G_k(y,\eta)-G_k(e_1,\eta))|\eta |^2(\mathfrak{h}_k(\delta_k\eta)e^{v_k(y)}-\mathfrak{h}_k(\delta_ke_1)e^{V_k})d\eta+O(\delta_k^2)$$
where the last term $O(\delta_k^2)$ comes from the oscillation of $w_k$ on $\partial \Omega_k$. Then we have
\begin{align*}
    w_k(y)
    &=-\frac 1{2\pi}\int_{\Omega_k}\log \frac{|y-\eta |}{ |e_1-\eta |}|\eta |^2(\mathfrak{h}_k(\delta_k \eta)-\mathfrak{h}_k(\delta_k e_1)e^{V_k})+O(\delta_k^2)\\
    &=-4\log \frac{|y-Q_2^k|}{|e_1-Q_2^k|}+4\log \frac{|y-e^{i\pi}|}{2}+O(\delta_k^2\mu_k).
\end{align*}
By $|Q_2^k-e^{i\pi}|=O(\delta_k^2)$ we see that $w_k(y)=O(\delta_k^2)$ on $|y-e^{i\pi}|=\tau$.

The standard point-wise estimate for singular equation ( see \cite{zhangcmp,gluck} ) gives 
\begin{align*}
&v_k(Q_2^k+\epsilon_kz)+2\log \epsilon_k\\
=&\log \frac{e^{\mu_k}}{(1+\frac{e^{\mu_k}}{8\mathfrak{h}_k(\delta_kQ_k)}|z|^2)^2}+\phi_1^k+C\delta_k^2\Delta (\log \mathfrak{h}_k)(0)(\log (1+|z|))^2, \quad |z|\sim \epsilon_k^{-1}. 
\end{align*}

\begin{align*}
&V_k(e^{i\pi}+\epsilon_kz)+2\log \epsilon_k \\
=&\log \frac{e^{\mu_k}}{(1+\frac{e^{\mu_k}}{8\mathfrak{h}_k(\delta_ke_1)}|z|^2)^2}+\phi_2^k+O(\epsilon_k^2(\log \epsilon_k)^2),\quad |z|\sim \epsilon_k^{-1}. 
\end{align*}

Thus 
$$w_k(Q_2^k+\epsilon_kz)=O(\epsilon_k^2(\log \epsilon_k)^2)+\phi_1^k-\phi_2^k+C\Delta (\log \mathfrak{h}_k)(0)\delta_k^2(\log (1+|z|)^2), $$
for $ |z|\sim \epsilon_k^{-1}. $
Taking the average around the origin, the spherical averages of the two harmonic functions are zero and $O(\delta_k^2)$ respectively, since they take zero at the origin and a point  at most $O(\delta_k^2)$ from the origin. So the spherical average of $w_k$ is comparable to 
$$\Delta (\log \mathfrak{h}_k)(0)\delta_k^2 (\log \epsilon_k)^2$$ for $|z|\sim \epsilon_k^{-1}$. Thus we know $\Delta (\log \mathfrak{h}_k)(0)=o(1)$ because $w_k=O(\delta_k^2\mu_k)$ in this region, 
Theorem \ref{main-thm-1} is established for all the cases. $\Box$

\section{Singular Mean field equation}

In this section we prove Theorem \ref{global-estimate}. First it is well known that if $p$ is a blowup point that has a non-quantized singular source ($\alpha_p=0$ or $\alpha_p\not \in \mathbb N$), the profile of the bubbling solutions around $p$ is a simple blwoup (see \cite{zhangcmp,zhangccm}).  So we only need to focus on the case that $\alpha_p\in \mathbb N$. Let $G(\cdot,\cdot)$ be the Green's function corresponding to $-\Delta_g$: 
$$-\Delta_y G(p,y)=\delta_p-1, \quad \int_M G(p,y)dV_g=0. $$
By setting 
$$G_1(y)=4\pi\sum_{t=1}^M\alpha_tG(p_t,y), $$
we have 
$$-\Delta G_1=4\pi \sum_t \alpha_t (\delta_{p_t}-1). $$
Then the function $v_k=u_k+G_1$ satisfies
$$\Delta_g v_k+\rho^k(\frac{h e^{v_k}e^{-G_1}}{\int_M he^{u_k}}-1)=0. $$
Let $p_1$ be a quantized singular source, which means $\alpha_{p_1}\in 4\pi \mathbb N$, in the neighborhood of $p_1$ we have
$$\Delta_g v_k+|y-p_1|^{2\alpha_{p_1}}H_ke^{v_k}=\rho_k $$
where
$$H_k=\frac{-\rho_khe^{4\pi \alpha_1\gamma(p_1,y)-4\pi\sum_{t\neq 1}\alpha_t G(p_t,y)}}{\int_M h e^{u_k}},$$
where $\gamma$ is the regular part of the Green's function. 
In local coordinates around $p_1$, the equation can be written as
$$\Delta v_k+|x|^{2\alpha_{p_1}}He^{\phi}e^{v_k}=\rho e^{\phi}. $$
where $\phi(0)=|\nabla \phi(0)|=0$ and $\Delta \phi(0)=-2K(p_1)$. 
Finally we use $f$ to remove the right hand side:
$$\Delta f=\rho^k e^{\phi}, \quad f(0)=0, \quad f=\mbox{constant on } \partial B_{\tau}$$
for $\tau>0$ small.  When we consider $v^k-f$ as the blowup solutions, we have 
$$\Delta (v_k-f)+|y|^{2\alpha_{p_1}}H_ke^fe^{v_k-f}=0. $$
It is a standard result that $H_k$ is uniformly bounded above and below.
From the defintion of $H_k$ we have
$$\Delta H_k(0)=\Delta h(p_1)-\sum_{t=1}^M4\pi \alpha_t. $$
Using Theorem \ref{main-thm-1} we would have
$$\Delta \log H_k(0)+\Delta \phi+\Delta f=o(1). $$
 if non-simple blowup happens at $p_1$, which is 
\begin{equation}\label{add-2}
\Delta \log h(p_1)-2K(p_1)-4\pi\sum_{t=1}^{M}\alpha_t+\rho^k=o(1). 
\end{equation}
Since $\rho^k\to \rho\in \Lambda$, we see from (\ref{curvature-a}) that (\ref{add-2}) cannot hold. Theorem \ref{global-estimate} is established. $\Box$
\section{Appendix}
In this section we prove (\ref{late-2-r}). Here we recall that $v_k$ is close to $V_{s,k}$ near $Q_s^k$ (see \ref{Q-bad})). That is why (\ref{late-1}) is used here.
The terms of $\phi_1$ and $\phi_3$ lead to $o(\epsilon_k)$, the integration involving $\phi_2$ cancels with the second term of (\ref{late-2-r}). The computation of $\phi_2$ is based on this equation:
$$\int_{\mathbb R^2}\frac{\frac{\mathfrak{h}_k(\delta_k Q_l^k)}4\sigma_k\nabla \mathfrak{h}_k(\delta_k Q_l^k)(Q_l^k-Q_s^k)|z|^2}{(1+\frac{\mathfrak{h}_k(\delta_k Q_l^k)}8|z|^2)^3}dz
=8\pi \sigma_k\nabla (\log \mathfrak{h}_k)(\delta_k Q_l^k)(Q_l^k-Q_s^k), $$
and by (\ref{delta-small-1})
\begin{equation}\label{subtle-1}
\nabla \log \mathfrak{h}_k(\delta_kQ_l^k)-\nabla \log \mathfrak{h}_k(\delta_k Q_s^k)=O(\delta_k)=o(\epsilon_k). 
\end{equation}

The integration involving $\phi_4$ provides the leading term. More detailed information is the following:
First for a global solution
$$V_{\mu,p}=\log \frac{e^{\mu}}{(1+\frac{e^{\mu}}{\lambda}|z^{N+1}-p|^2)^2}$$ of
$$\Delta V_{\mu,p}+\frac{8(N+1)^2}{\lambda}|z|^{2N}e^{V_{\mu,p}}=0,\quad \mbox{in }\quad \mathbb R^2, $$
by differentiation with respect to $\mu$ we have
$$\Delta(\partial_{\mu}V_{\mu,p})+\frac{8(N+1)^2}{\lambda}|z|^{2N}e^{V_{\mu,p}}\partial_{\mu}V_{\mu,p}=0,\quad \mbox{in}\quad \mathbb R^2. $$
By the expression of $V_{\mu,p}$ we see that 
$$\partial_r\bigg (\partial_{\mu}V_{\mu,p}\bigg )(x)=O(|x|^{-2N-3}).$$ 
Thus we have
\begin{equation}\label{inte-eq-1}
\int_{\mathbb R^2}\partial_{\mu}V_{\mu,p}|z|^{2N}e^{V_{\mu,p}}=\int_{\mathbb R^2}\frac{(1-\frac{e^{\mu}}{\lambda}|z^{N+1}-P|^2)|z|^{2N}}{(1+\frac{e^{\mu}}{\lambda}|z^{N+1}-P|^2)^3}dz=0.
\end{equation}

From $V_{\mu,p}$  we also have
$$\int_{\mathbb R^2}\partial_{P}V_{\mu,p}|y|^{2N}e^{V_{\mu,p}}=\int_{\mathbb R^2}\partial_{\bar P}V_{\mu,p}|y|^{2N}e^{V_{\mu,p}}=0, $$
which gives
\begin{equation}\label{inte-eq-2}
\int_{\mathbb R^2}\frac{\frac{e^{\mu}}{\lambda}(\bar z^{N+1}-\bar P)|z|^{2N}}{(1+\frac{e^{\mu}}{\lambda}|z^{N+1}-P|^2)^3}=\int_{\mathbb R^2}\frac{\frac{e^{\mu}}{\lambda}( z^{N+1}- P)|z|^{2N}}{(1+\frac{e^{\mu}}{\lambda}|z^{N+1}-P|^2)^3}=0.
\end{equation}

Now we need more precise expressions of $\phi_1$, $\phi_3$ and $B$:
\begin{align*}
&\phi_1=(\mu_s^k-\mu_l^k)(1-\frac{(N+1)^2}{D_s^k}|z+\frac{N}2\epsilon_k z^2e^{-i\beta_s}|^2/B, \\
&\phi_3=\frac{4(N+1)}{D_s^k B}Re((z+\frac{N}2\epsilon_ke^{-i\beta_s}z^2))(\frac{\bar p_s^k-\bar p_l^k}{\epsilon_k}e^{-i\beta_s}))
\\
&B=1+\frac{(N+1)^2}{D_s^k}|z+\frac N2z^2e^{-i\beta_s}\epsilon_k|^2,
\end{align*}
From here we use scaling and cancellation to have
$$\int_{B(0,\tau\epsilon_k^{-1})}\frac{\phi_1}{M_k}B^{-2}=o(\epsilon_k),\quad
\int_{B(0,\tau\epsilon_k^{-1})}\frac{\phi_3}{M_k}B^{-2}=o(\epsilon_k).$$
Thus (\ref{late-2-r}) holds.

\end{document}